\documentclass[12pt,a4paper]{article}
\usepackage[latin1]{inputenc}
\usepackage[T1]{fontenc}
\usepackage[english]{babel}
\usepackage{amssymb,amsmath}
\usepackage{bbm}

\newcommand{\thmheadercommand}[1]{\textbf{\scshape{}#1}}

\def\d{{\mathrm{d}}}
\def\N{{\mathbb{N}}}
\def\Z{{\mathbb{Z}}}

\def\R{{\mathbb{R}}}

\renewcommand{\geq}{\geqslant}
\renewcommand{\leq}{\leqslant}

\def\eps{\varepsilon}
\renewcommand{\epsilon}{\varepsilon}
\renewcommand{\phi}{\varphi}

\DeclareMathOperator{\diam}{diam}
\DeclareMathOperator{\vol}{vol}
\DeclareMathOperator{\Ric}{Ric}

\DeclareMathOperator{\Ent}{Ent}
\DeclareMathOperator{\Var}{Var}

\DeclareMathOperator{\Tr}{Tr}

\DeclareMathOperator{\Id}{Id}

\newcommand{\abs}[1]{\left|\mskip1mu#1\right|}
\newcommand{\norm}[1]{\left\|#1\right\|}

\newcounter{prop}
\newcounter{defi}
\newcounter{thm}
\newcounter{lem}

\newenvironment{dem}[1][]{\noindent{\thmheadercommand{Proof#1}}\,\,--\,\,}{$\square$\medskip}

\newenvironment{enonce}[1]{\medskip\noindent{\thmheadercommand{#1}}\,\,--\,\,\begin{slshape}}{\end{slshape}\medskip}

\newenvironment{enonce2}[1]{\medskip\noindent{\thmheadercommand{#1}}\,\,--\,\,}{\medskip}

\newenvironment{defi}[1][]{\refstepcounter{prop}
\begin{enonce}{Definition \theprop{}#1}}{\end{enonce}}
\newenvironment{prop}[1][]{\refstepcounter{prop}
\begin{enonce}{Proposition \theprop{}#1}}{\end{enonce}}
\newenvironment{thm}[1][]{\refstepcounter{prop}
\begin{enonce}{Theorem \theprop{}#1}}{\end{enonce}}
\newenvironment{lem}[1][]{\refstepcounter{prop}
\begin{enonce}{Lemma \theprop{}#1}}{\end{enonce}}
\newenvironment{cor}[1][]{\refstepcounter{prop}
\begin{enonce}{Corollary \theprop{}#1}}{\end{enonce}}
\newenvironment{ex}[1][]{\refstepcounter{prop}
\begin{enonce}{Example \theprop{}#1}}{\end{enonce}}

\newenvironment{rem}[1][]{\refstepcounter{prop}
\begin{enonce2}{Remark \theprop{}#1}}{\end{enonce2}}

\newenvironment{defi*}[1][]{
\begin{enonce}{Definition#1}}{\end{enonce}}
\newenvironment{prop*}[1][]{
\begin{enonce}{Proposition#1}}{\end{enonce}}
\newenvironment{thm*}[1][]{
\begin{enonce}{Theorem#1}}{\end{enonce}}
\newenvironment{lem*}[1][]{
\begin{enonce}{Lemma#1}}{\end{enonce}}
\newenvironment{cor*}[1][]{
\begin{enonce}{Corollary#1}}{\end{enonce}}
\newenvironment{ex*}[1][]{
\begin{enonce}{Example#1}}{\end{enonce}}
\newenvironment{exo*}[1][]{
\begin{enonce2}{Exercise#1}}{\end{enonce2}}
\newenvironment{rem*}[1][]{
\begin{enonce2}{Remark#1}}{\end{enonce2}}

\usepackage{graphicx}
\usepackage[dvips]{color}
\usepackage{url}

\title{Ricci curvature of Markov chains on metric spaces}

\author{Yann Ollivier}

\date{}

\newcommand{\E}{\mathbb{E}}
\DeclareMathOperator{\Supp}{Supp}
\newcommand{\Td}{\mathcal{T}_1}
\newcommand{\Tdd}{\mathcal{T}_2}
\newcommand{\Pm}{\mathcal{P}}
\renewcommand{\d}{\mathrm{d}\hspace{-0.02em}}
\newcommand{\deltat}{\hspace{0.05em}\delta\hspace{-.06em}t\hspace{0.05em}}
\newcommand{\e}{\mathrm{e}}

\begin{document}

\maketitle

\begin{abstract}
We define the Ricci curvature of Markov chains on metric spaces as a
local contraction coefficient of the random walk acting on the space of
probability measures equipped with a Wasserstein transportation distance.
For Brownian motion on a Riemannian manifold this gives back the value of
Ricci curvature of a tangent vector. Examples of positively curved spaces
for this definition include the discrete cube and discrete
versions of the Ornstein--Uhlenbeck process. Moreover this generalization
is consistent with the Bakry--Émery Ricci curvature for Brownian motion
with a drift on a Riemannian manifold.

Positive Ricci curvature is shown to imply a spectral gap, a
Lévy--Gromov-like Gaussian concentration theorem and a kind of
modified logarithmic Sobolev inequality. The bounds obtained are sharp in
several interesting examples.
\end{abstract}

\section*{Introduction}

There are numerous generalizations of the notion of a metric space with
\emph{negative sectional curvature}: manifolds with negative sectional
curvature, $\text{CAT}(0)$ and $\text{CAT}(-1)$ spaces or
$\delta$-hyperbolic spaces are widely used in various branches of
mathematics and give rise to numerous theorems. For \emph{positive}
curvature in Riemannian geometry, the right concept seems to be a lower
bound on \emph{Ricci} curvature (which is weaker than a lower bound on
sectional curvature). The most basic result in this direction is the
Bonnet--Myers theorem bounding the diameter of the space in function of
the Ricci curvature, but let us mention Lichnerowicz' theorem for the
spectral gap of the Laplacian (Theorem~181 in~\cite{Ber03}), the
Lévy--Gromov theorem for isoperimetric inequalities and concentration of
measure~\cite{Gro86}, or Gromov's theorem on precompactness of the space of
manifolds with given dimension, upper bound on the diameter and
lower bound on the Ricci curvature.

We refer to the nice survey~\cite{Lott} for a discussion of the geometric
interest of lower bounds on Ricci curvature, with further references, and
the need for a generalized notion of positive Ricci curvature for metric
spaces (often equipped with a measure).

There have been several generalizations of the notion of Ricci curvature.
First, the study by Bakry and Émery~\cite{BE85} of hypercontractivity of
diffusion processes led them to show that, when considering the Brownian
motion on a manifold with an additional drift given by a tangent vector
field $F$, the quantity $\Ric-2\nabla^{\text{sym}}F$ plays the role of a
Ricci curvature for the process, as far as functional inequalities are
concerned. The main example is the Ornstein--Uhlenbeck process on $\R^N$,
whose invariant distribution is Gaussian, and which is positively curved
in this sense.

Later, simultaneously, Sturm~\cite{Stu06}, Lott and Villani~\cite{LV},
and Ohta~\cite{Oht} used ideas from optimal transportation theory to
define a notion of lower bound on the Ricci curvature for length spaces
equipped with a measure.  Their definition keeps a lot of the properties
traditionally associated with positive Ricci curvature, and is compatible
with the Bakry--Émery extension. However, it has two main drawbacks.
First, it is infinitesimal, and in particular is meaningless for a graph.
Second, the definition is rather involved and difficult to check on
concrete examples. The main class of spaces for which this definition is
interesting are Gromov--Hausdorff limits of manifolds of a given
dimension.

Here we propose a definition of Ricci curvature for metric spaces
equipped with a Markov chain or a diffusion process (which for a
Riemannian manifold will typically be Brownian motion), which is
hopefully simpler to check on examples. The definition is again based on
optimal transportation, but in a less infinitesimal way, and can be used
to define a notion of ``curvature at a given scale'' for a metric space.
As a consequence, we can test it in discrete spaces such as graphs.  Such
an example is the discrete cube $\{0,1\}^N$, which from the point of view
of concentration of measure behaves very much like the sphere $S^N$, and
is thus expected to somehow have positive curvature.

Our definition, when applied to a Riemannian manifold equipped with the
Brownian motion, gives back the usual value of the Ricci curvature of a
tangent vector. It is consistent with the Bakry--Émery extension, and
provides a visual explanation for the contribution $-\nabla^{\text{sym}}
F$ of the drift $F$.  We are able to prove generalizations of the
Bonnet--Myers theorem, of the Lichnerowicz spectral gap theorem and of
the Lévy--Gromov isoperimetry theorem, as well as a kind of modified logarithmic
Sobolev inequality, although with some (bounded) loss in the constants.
As a by-product, we get a new proof for Gaussian concentration and the
logarithmic Sobolev inequality in the Lévy--Gromov or Bakry--Émery
context (though the constants are not sharp).

\paragraph{Related work.} After having written a first version of this
text, we learned that related ideas appear in several recent papers.
Joulin~\cite{Jou} uses contraction of the Lipschitz constant (under the
name ``Wasserstein curvature'') to get a Poisson-type concentration
result for continuous-time Markov chains on a countable space, at least
in the bounded, one-dimensional case. Oliveira~\cite{Oli} proves that
Kac's random walk on $\mathrm{SO}(n)$ has positive Ricci curvature in our
sense, which allows to improve mixing time estimates significantly.
Djellout, Guillin and Wu~\cite{DGW04} use contraction of Lipschitz
constants and transportation distances (without the link with Ricci
curvature) in the context of dependent sequences of random variables, to
get Gaussian concentration results.  The link with the spectral gap
appears in~\cite{Sam} (p.~94) for the particular case of graphs, and is
present in the works of Chen (e.g.~\cite{CW97,Che98}).

From the discrete Markov chain point of view, the techniques presented
here are just a metric version of the usual coupling method. Namely,
Ricci curvature can be seen as a refined version of Dobrushin's ergodic
coefficient (see \cite{Dob56}, or e.g.\ section~6.7.1 in~\cite{Bre99})
using the metric structure on the underlying space.

From the Riemannian point of view, our approach boils down to contraction
of the Lipschitz norm by the heat equation, which is one of the results
of Bakry and Émery (\cite{BE84,BE85}, see also~\cite{LSI}
and~\cite{RS05}). This latter property was suggested in~\cite{RS05} as a
possible definition of a lower bound on Ricci curvature for diffusion
operators in general spaces, though it does not provide an explicit value
for Ricci curvature at a given point.

\paragraph{Acknowledgements.} I would like to thank Vincent Beffara,
Fabrice Debbasch, Alessio Figalli, Pierre Pansu, Bruno Sévennec, Romain
Tessera and Cédric Villani for numerous inspiring conversations about
coarse geometry and Ricci curvature,  as well as Djalil Chafaï, Aldéric
Joulin, Shin-ichi Ohta and Roberto Oliveira for useful remarks on the
manuscript and bibliographical references. Special thanks to Pierre Py
for the two points $x$ and $y$.

\paragraph{Notation.} In the paper, we use the symbol $\approx$ to denote
equality up to a multiplicative universal constant (typically $2$ or
$4$); the symbol $\sim$ denotes usual asymptotic equivalence. The word
``distribution'' is used as a synonym for ``probability measure''.

\section{Definitions and statements}

\subsection{Ricci curvature}

A common framework for generalizations of Ricci curvature is that of
metric measure spaces~\cite{Stu06,LV}. However, most measures appear as
the invariant distribution of some process (e.g.\ Brownian motion on a Riemannian
manifold), and it is more convenient and more general to start with a
process in a metric space, as is the case in Bakry--Émery theory. See also
Remark~\ref{rem:randomwalk} below.

Here for simplicity we will mainly consider the case of a discrete-time
process. Similar definitions and results can be given for continuous
time (see e.g.\ Section~\ref{sec:continuoustime}).

\begin{defi}
\label{def:space}
Let $(X,d)$ be a Polish metric space, equipped with its Borel
$\sigma$-algebra.

A \emph{random walk $m$ on
$X$} is a family of probability measures $m_x(\cdot)$ on $X$ for each
$x\in X$, satisfying the following two technical assumptions: $(i)$ the
measure $m_x$ depends measurably on the point $x\in X$; $(ii)$ each
measure $m_x$ has finite first moment, i.e.\ for some (hence any) $o\in
X$ one has $\int d(o,y)\,\d m_x(y)<\infty$.
\end{defi}

This defines a Markov chain whose transition probability from $x$ to $y$
in $n$ steps is
\[
\d m_x^{\ast n}(y):=\int_{z\in X} \d m_x^{\ast (n-1)}(z)\,\d m_z(y)
\]
where of course $m_x^{\ast 1}:=m_x$.

Recall that a measure $\nu$ on $X$ is \emph{invariant} for this random
walk if $\d\nu(x)=\int_y \d\nu(y) \d m_y(x)$. It is \emph{reversible} if
moreover, the detailed balance condition $\d\nu(x)\d m_x(y)=\d\nu(y) \d
m_y(x)$ holds.

\bigskip

This allows to define a notion of curvature as follows. Consider two very
close points $x,y$ in a Riemannian manifold, defining a tangent vector
$(xy)$. Let $w$ be another tangent vector at $x$; let $w'$ be the tangent
vector at $y$ obtained by parallel transport of $w$ from $x$ to $y$. Now
if we follow the two geodesics issuing from $x,w$ and $y,w'$, in positive
curvature the geodesics will get closer, and will part away in negative
curvature.  Ricci curvature along $(xy)$ is this phenomenon, averaged on
all directions $w$ at $x$.

So in the general case, we will measure whether following the random walk
issuing from two nearby points $x,y$ results in points that are closer
than $x,y$ were, in which case Ricci curvature will be positive, or
further apart, in which case Ricci curvature will be negative.  This is
made precise by the use of transportation distances between
probability measures. We refer to~\cite{Vil03} for an introduction to
this topic.

\begin{defi}
Let $(X,d)$ be a metric space and let $\nu_1, \nu_2$ be two probability
measures on $X$. The \emph{$L^1$ transportation distance} between $\nu_1$
and $\nu_2$ is
\[
\Td(\nu_1,\nu_2):=\inf_{\xi\in \Pi(\nu_1,\nu_2)} \int_{(x,y)\in X\times X}
d(x,y)\, \d\xi(x,y)
\]
where $\Pi(\nu_1,\nu_2)$ is the set of
measures on $X\times X$ projecting to $\nu_1$ and $\nu_2$.
\end{defi}

Intuitively, $\d\xi(x,y)$ represents the mass that is sent from $x$ to
$y$, hence the constraint on the projections of $\xi$, ensuring that the
initial measure is $\nu_1$ and the final measure is $\nu_2$.

The infimum is actually attained (Theorem~1.3 in~\cite{Vil03}), but the
optimal coupling is generally not unique. In what follows, it is enough
to chose one such coupling.

\begin{defi}
Let $(X,d)$ be a metric space with a random walk $m$. Let $x,y\in X$ be
two distinct points. The \emph{Ricci curvature of $(X,d,m)$ in the
direction $(x,y)$} is
\[
\kappa(x,y):=1-\frac{\Td(m_x,m_y)}{d(x,y)}
\]
\end{defi}

\begin{center}
\begin{picture}(0,0)%
\includegraphics{defricci.pstex}%
\end{picture}%
\setlength{\unitlength}{4144sp}%
\begingroup\makeatletter\ifx\SetFigFont\undefined%
\gdef\SetFigFont#1#2#3#4#5{%
  \reset@font\fontsize{#1}{#2pt}%
  \fontfamily{#3}\fontseries{#4}\fontshape{#5}%
  \selectfont}%
\fi\endgroup%
\begin{picture}(4103,3197)(952,-3044)
\put(1261,-2986){\makebox(0,0)[lb]{\smash{\SetFigFont{12}{14.4}{\rmdefault}{\mddefault}{\updefault}{\color[rgb]{0,0,0}$x$}%
}}}
\put(2431,-2671){\makebox(0,0)[lb]{\smash{\SetFigFont{12}{14.4}{\rmdefault}{\mddefault}{\updefault}{\color[rgb]{0,0,0}$d(x,y)$}%
}}}
\put(3961,-2986){\makebox(0,0)[lb]{\smash{\SetFigFont{12}{14.4}{\rmdefault}{\mddefault}{\updefault}{\color[rgb]{0,0,0}$y$}%
}}}
\put(1261,-1996){\makebox(0,0)[lb]{\smash{\SetFigFont{12}{14.4}{\rmdefault}{\mddefault}{\updefault}{\color[rgb]{0,0,0}$m_x$}%
}}}
\put(4096,-1951){\makebox(0,0)[lb]{\smash{\SetFigFont{12}{14.4}{\rmdefault}{\mddefault}{\updefault}{\color[rgb]{0,0,0}$m_y$}%
}}}
\put(2431,-2041){\makebox(0,0)[lb]{\smash{\SetFigFont{12}{14.4}{\rmdefault}{\mddefault}{\updefault}{\color[rgb]{0,0,0}on average}%
}}}
\put(2341,-1816){\makebox(0,0)[lb]{\smash{\SetFigFont{12}{14.4}{\rmdefault}{\mddefault}{\updefault}{\color[rgb]{0,0,0}$(1-\kappa)d(x,y)$}%
}}}
\end{picture}

\end{center}

When $(X,d)$ is a Riemannian manifold, if the random walk consists in
randomly jumping in a ball of radius $\eps$ around $x$, for small $\eps$
and close enough $x,y$ this definition captures the Ricci
curvature in the direction $xy$ (up to some factor depending on $\eps$).

We will see below (Proposition~\ref{prop:epsgeod}) that in geodesic spaces, it is enough to know
$\kappa(x,y)$ for close points $x,y$.

If a continuous-time Markov kernel is given,
one can also define a continuous-time version of the Ricci curvature by setting
\[
\kappa(x,y):=-\,\frac{\d}{\d t} \,\frac{\Td(m^t_x,m^t_y)}{d(x,y)}
\]
when this derivative exists, but for simplicity we will
mainly work with the discrete-time version here. Indeed, for
continuous-time Markov chains, existence of the process is already a
non-trivial issue. We will sometimes use our results on concrete
continuous-time examples (e.g.~$M/M/\infty$ queues in
section~\ref{sec:continuoustime}), but only when they appear as
an obvious limit of a discrete-time approximation.

One could use the $L^p$ transportation distance instead of the $L^1$ one
in the definition; however, though this will result in stronger
assumptions, I did not find any theorem where this would be necessary.

\begin{enonce2}{Notation }
By analogy with the Riemannian case, when computing the transportation
distance between measures $m_x$ and $m_y$, we will think of $X\times X$
equipped with the coupling measure as a tangent space, and for
$z\in X\times X$ we will write $x+z$ and $y+z$ for the two projections to
$X$. So in this notation we have
\[
\kappa(x,y)=-\,\frac{1}{d(x,y)}\int (d(x+z,y+z)-d(x,y))\,\d z
\]
where implicitly $\d z$ is the optimal coupling between $m_x$ and $m_y$.
\end{enonce2}

\subsection{Examples}

\begin{ex}[ ($\Z^N$ and $\R^N$)]
Let $m$ be the simple random walk on the graph of the grid $\Z^N$
equipped with its graph metric. Then for any two points $x,y\in \Z^d$,
the Ricci curvature along $(xy)$ is $0$.
\end{ex}

Indeed, we can transport the measure $m_x$ around $x$ to the measure
$m_y$ by a translation of vector $y-x$ (and this is optimal), so that the
distance between $m_x$ and $m_y$ is exactly that between $x$ and $y$.

This example generalizes to the case of $\Z^n$ or $\R^N$ equipped with
any translation-invariant norm and any random walk given by a
translation-invariant transition kernel (consistently with~\cite{LV}).
For example, the triangular tiling of the plane has $0$ curvature.

\begin{rem}[ (Random walk at scale $\eps$)]
\label{rem:randomwalk}
It is easy to construct random walks on metric measure spaces.
If $(X,d,\mu)$ is a metric measure space (for example with $\mu$ the
Hausdorff measure) and $\eps>0$, the \emph{random walk at scale $\eps$}
consists in, starting at a point $x$, randomly jumping in the ball of
radius $\eps$ around $x$, with probability density proportional to $\mu$;
namely $\d m_x(y):=\d\mu(y)/\mu(B(x,\eps))$ if $d(x,y)\leq \eps$ (one can
also use other functions of the distance, such as a Gaussian kernel).
This allows to consider the Ricci curvature associated with this random
walk.
\end{rem}

This is what we do now on Riemannian manifolds to get back the
usual Ricci curvature (up to some normalization constants), hence the
terminology.

\begin{prop}
\label{prop:riemdist}
Let $(X,d)$ be a smooth complete Riemannian manifold. Let $v,w$ be
unit tangent
vectors at $x\in X$. Let $\eps,\delta>0$. Let $y=\exp_x \delta v$ and let
$w'$ be the tangent vector at $y$ obtained by parallel transport of $w$ along the
geodesic $\exp_x tv$. Then
\[
d(\exp_x \eps w, \exp_y \eps w')= \delta
\left(1-\frac{\eps^2}{2}\left(K(v,w)+O(\delta+\eps)\right)\right)
\]
as $(\eps,\delta)\to 0$. Here $K(v,w)$ is the sectional curvature in the
tangent plane $(v,w)$.
\end{prop}

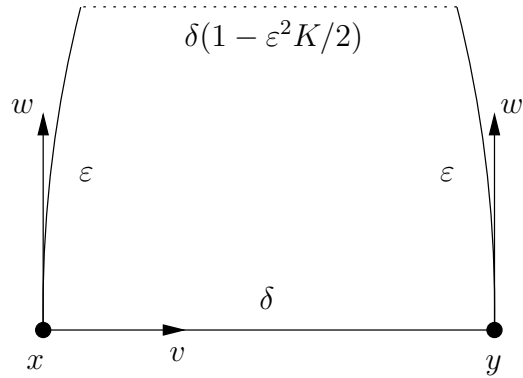
\begin{center}
\begin{picture}(0,0)%
\includegraphics{Ksect.pstex}%
\end{picture}%
\setlength{\unitlength}{4144sp}%
\begingroup\makeatletter\ifx\SetFigFont\undefined%
\gdef\SetFigFont#1#2#3#4#5{%
  \reset@font\fontsize{#1}{#2pt}%
  \fontfamily{#3}\fontseries{#4}\fontshape{#5}%
  \selectfont}%
\fi\endgroup%
\begin{picture}(2933,2230)(1171,-3044)
\put(1261,-2986){\makebox(0,0)[lb]{\smash{\SetFigFont{12}{14.4}{\rmdefault}{\mddefault}{\updefault}{\color[rgb]{0,0,0}$x$}%
}}}
\put(4096,-1456){\makebox(0,0)[lb]{\smash{\SetFigFont{12}{14.4}{\rmdefault}{\mddefault}{\updefault}{\color[rgb]{0,0,0}$w'$}%
}}}
\put(1171,-1456){\makebox(0,0)[lb]{\smash{\SetFigFont{12}{14.4}{\rmdefault}{\mddefault}{\updefault}{\color[rgb]{0,0,0}$w$}%
}}}
\put(1576,-1861){\makebox(0,0)[lb]{\smash{\SetFigFont{12}{14.4}{\rmdefault}{\mddefault}{\updefault}{\color[rgb]{0,0,0}$\eps$}%
}}}
\put(3736,-1861){\makebox(0,0)[lb]{\smash{\SetFigFont{12}{14.4}{\rmdefault}{\mddefault}{\updefault}{\color[rgb]{0,0,0}$\eps$}%
}}}
\put(2206,-1051){\makebox(0,0)[lb]{\smash{\SetFigFont{12}{14.4}{\rmdefault}{\mddefault}{\updefault}{\color[rgb]{0,0,0}$\delta (1-\eps^2 K/2)$}%
}}}
\put(4006,-2986){\makebox(0,0)[lb]{\smash{\SetFigFont{12}{14.4}{\rmdefault}{\mddefault}{\updefault}{\color[rgb]{0,0,0}$y$}%
}}}
\put(2116,-2941){\makebox(0,0)[lb]{\smash{\SetFigFont{12}{14.4}{\rmdefault}{\mddefault}{\updefault}{\color[rgb]{0,0,0}$v$}%
}}}
\put(2656,-2626){\makebox(0,0)[lb]{\smash{\SetFigFont{12}{14.4}{\rmdefault}{\mddefault}{\updefault}{\color[rgb]{0,0,0}$\delta$}%
}}}
\end{picture}

\end{center}

\begin{ex}[ (Riemannian manifold)]
\label{ex:riem}
Let $(X,d)$ be a smooth complete $N$-dimensional Riemannian manifold. For some $\eps>0$, let the
Markov chain $m^\eps$ be defined by
\[
\d m^\eps_x(y):=\frac{1}{\vol(B(x,\eps))} \,\d\vol(y)
\]
if $y\in B(x,\eps)$, and $0$ otherwise.

Let $x\in X$ and let $v$ be a unit tangent vector at $x$. Let $y$ be a
point on the geodesic issuing from $v$, with $d(x,y)$ small enough.
Then
\[
\kappa(x,y)=\frac{\eps^2}{2(N+2)}\left(\Ric(v,v)+O(\eps)+O(d(x,y))\right)
\]
\end{ex}

\begin{dem}
This is essentially the same as Theorem~1.5 (condition $(xii)$)
in~\cite{RS05}, except that
therein, the infimum of Ricci curvature is used instead of its value
along a tangent vector.
The proof is postponed to Section~\ref{sec:riemproof}.
Basically, the value of $\kappa(x,y)$ is obtained by averaging the
proposition above for $w$ in the unit ball of the tangent space at $x$, which provides an
upper bound for $\kappa$. The lower bound requires use of the dual
characterization of transportation distance (Theorem~1.14
in~\cite{Vil03}).
\end{dem}

\begin{ex}[ (Discrete cube)]
\label{ex:cube}
Let $X=\{0,1\}^N$ be the discrete cube equipped with the Hamming metric
(each edge is of length $1$). Let $m$ be the lazy random walk on the
graph $X$,
i.e.\ $m_x(x)=1/2$ and $m_x(y)=1/2N$ if $y$ is a neighbor of $x$.

Let $x,y\in X$ be neighbors. Then $\kappa(x,y)=1/N$.
\end{ex}

This examples generalizes to arbitrary binomial distributions (see
Section~\ref{sec:poisson}).

Here laziness is necessary to avoid parity problems: If no laziness is
introduced, points at odd distance never meet under the random walk; in
this case one must consider Ricci curvature for points at even
distance only.

Actually, since the discrete cube is a $1$-geodesic space,
one has $\kappa(x,y)\geq 1/N$ for any pair $x,y\in X$, not only neighbors
(see Proposition~\ref{prop:epsgeod}).

\begin{dem}
We can suppose that $x=00\ldots0$ and $y=10\ldots0$. For $z\in X$ and
$1\leq i\leq N$, let us denote by $z^i$ the neighbor of $z$ in which the
$i$-th bit is switched. An optimal coupling between $m_x$ and $m_y$ is as
follows: For $i\geq 2$, move $x^i$ to $y^i$ (both have mass $1/2N$ under
$m_x$ and $m_y$ respectively). Now $m_x(x)=1/2$ and $m_y(x)=1/2N$, and
likewise for $y$.  To transport $m_x$ to $m_y$, it is enough to move a
mass $1/2-1/2N$ from $x$ to $y$. All points are moved over a distance $1$
by this coupling, except for a mass $1/2N$ which remains at $x$ and a
mass $1/2N$ which remains at $y$, and so the Ricci curvature is at least
$1/N$.

Optimality of this coupling is obtained as follows: Consider the function
$f:X\to \{0,1\}$ which sends a point of $X$ to its first bit. This is a
$1$-Lipschitz function, with $f(x)=0$ and $f(y)=1$.
The expectations of $f$ under $m_x$ and $m_y$ are $1/2N$ and $1-1/2N$
respectively, so that $1-1/N$ is a lower bound on $\Td(m_x,m_y)$.

A very short but less visual proof can be obtained through the
$L^1$ tensorization property (Proposition~\ref{prop:tensorization}).
\end{dem}

\begin{ex}[ (Ornstein--Uhlenbeck process)]
\label{ex:OU}
Let $s\geq 0, \alpha>0$ and consider the Ornstein--Uhlenbeck process in
$\R^N$ given by the stochastic differential equation
\[
\d X_t=-\alpha X_t \,\d t + s\, \d B_t
\]
where $B_t$ is a standard $N$-dimensional Brownian motion. The invariant
distribution is Gaussian, of variance $s^2/2\alpha$.

Let $\deltat>0$ and let the random walk $m$ be the flow at time $\delta
t$ of
the process. Explicitly, $m_x$ is a Gaussian probability measure centered
at $\e^{-\alpha\deltat}x$, of variance
$s^2(1-\e^{-\alpha\deltat})/\alpha\sim s^2\deltat$ for small $\deltat$.

Then the Ricci curvature $\kappa(x,y)$ of this random walk is
$1-\e^{-\alpha \deltat}$, for any two $x,y\in \R^N$.
\end{ex}

\begin{dem}
The transportation distance between two Gaussian distributions with the
same variance is the distance between their centers, so that
$\kappa(x,y)=1-\frac{\abs{\e^{-\alpha\deltat}x-\e^{-\alpha\deltat}y}}{\abs{x-y}}$.
\end{dem}

\begin{ex}[ (Discrete Ornstein--Uhlenbeck)]
\label{ex:discrOU}
Let $X=\{-N,-N+1,\ldots,N-1,N\}$ and let $m$ be the random walk on $X$ given by
\[
m_k(k)=1/2,\qquad m_k(k+1)=1/4-k/4N,\qquad m_k(k-1)=1/4+k/4N
\]
which is a lazy random walk with linear drift towards $0$. The binomial
distribution $\frac1{2^{2N}} \binom{2N}{N+k}$ is reversible for this
random walk.

Then, for any two neighbors $x,y$ in $X$, one has $\kappa(x,y)=1/2N$.
\end{ex}

\begin{dem}
Exercise.
\end{dem}

\begin{ex}[ (Bakry--Émery)]
\label{ex:BE}
Let $X$ be an $N$-dimensional Riemannian manifold and $F$ be a tangent vector field. Consider the differential operator
\[
L:=\frac12 \Delta+F.\nabla
\]
associated with the stochastic differential equation
\[
\d x_t=F\,\d t+\d B_t
\]
where $B_t$ is the Brownian motion in $X$. The Ricci curvature (in the
Bakry--Émery sense) of this operator is $\frac12\Ric-\nabla^{\text{sym}}F$ where
$\nabla^{\text{sym}}F^{ij}=\frac12(\nabla^iF^j+\nabla^jF^i)$ is the
symmetrized of $\nabla F$.

Consider the Euler approximation scheme at time $\deltat$ for this
stochastic equation, which consists in following the flow of $F$ for a
time $\deltat$ and then randomly jumping in a ball of radius
$\sqrt{(N+2)\deltat}$.

Let $x\in X$ and let $v$ be a unit tangent vector at $x$. Let $y$ be a
point on the geodesic issuing from $v$, with $d(x,y)$ small enough.
Then
\[
\kappa(x,y)=\deltat\left(\frac12
\Ric(v,v)-\nabla^{\text{sym}}F(v,v)+O(d(x,y))+O(\sqrt{\deltat})\right)
\]
\end{ex}

\begin{dem}
First let us explain the normalization: Jumping in a ball of radius
$\eps$ generates a variance $\eps^2\frac{1}{N+2}$ in a given direction. On the other hand,
the $N$-dimensional Brownian motion has, by definition, a variance $\d t$
per unit of time $\d t$ in any given direction, so a proper discretization at time $\deltat$
requires jumping in a ball of radius $\sqrt{(N+2)\deltat}$.
Also, as noted in~\cite{BE85}, the generator of Brownian
motion is $\frac12 \Delta$ instead of $\Delta$, hence the
$\frac12$ factor for the Ricci part.

Now the discrete-time process begins by following the flow $F$ for some
time $\deltat$. Starting at points $x$ and $y$, using elementary
Euclidean geometry, it is easy to see that after this, the distance
between the endpoints behaves like
$d(x,y)(1+\deltat\,v.\nabla_{\!\!v}F+O(\deltat^2))$. Note that
$v.\nabla_{\!\!v}F=\nabla^{\text{sym}}F(v,v)$.

Now, just as in Example~\ref{ex:riem}, randomly jumping in a ball of radius
$\eps$ results in a gain of $d(x,y)\frac{\eps^2}{2(N+2)}\Ric(v,v)$ on
transportation distances. Here
$\eps^2=(N+2)\deltat$. So after
the two steps, the distance between the endpoints is
\[
d(x,y)\left(1-\frac{\deltat}{2}\Ric(v,v)+\deltat
\,\nabla^{\text{sym}}F(v,v)\right)
\]
as needed, up to higher-order terms.
\end{dem}

Maybe the reason for the additional $-\nabla^{\text{sym}} F$ in Ricci
curvature à la Bakry--Émery is made clearer in this context: it is simply
the quantity by which the flow of $X$ modifies distances between two
starting points.

It is clear on this example why reversibility is not fundamental in this
theory: the antisymmetric part of the force $F$ generates an
infinitesimal isometric displacement. Combining the Markov chain with an
isometry of the space has no effect whatsoever on our definition.

\begin{ex}[ (Multinomial distribution)]
\label{ex:mult}
Consider the set $X=\{(x_0,x_1,\ldots,x_d),\, x_i\in \N, \,\sum x_i=N\}$ viewed
as the configuration set of $N$ balls in $d+1$ boxes. Consider the
process which consists in taking a ball at random among the $N$ balls,
removing it from its box, and putting it back at random in one of the $d+1$
boxes. More precisely, the transition probability from $(x_0,\ldots,x_d)$
to $(x_0,\ldots,x_i-1,\ldots,x_j+1,\ldots,x_d)$ (with maybe $i=j$) is
$x_i/N(d+1)$.  The multinomial distribution $\frac{N!}{(d+1)^N\,\prod x_i!}$ is reversible
for this Markov chain.

Equip this configuration space with the metric $d((x_i),(x'_i)):=\frac12
\sum \abs{x_i-x'_i}$ which is the graph distance w.r.t.\ the moves above.
Then the Ricci curvature of the Markov chain is $1/N$.
\end{ex}

\begin{dem}
Exercise.
\end{dem}

\begin{ex}[ (Geometric distribution)]
\label{ex:geom}
Let the random walk on $\N$ be defined by the transition probabilities
$p_{n,n+1}=1/3$, $p_{n+1,n}=2/3$ and $p_{0,0}=2/3$. This random walk is
reversible with respect to the geometric measure $2^{-(n+1)}$. It is
easy to check that for $n\geq 1$ one has $\kappa_{n,n+1}=0$.
\end{ex}

\begin{dem}
The transition kernel is translation-invariant except at $0$.
\end{dem}

Section~\ref{sec:exp} contains more material about this latter example and how
non-negative Ricci curvature sometimes implies exponential concentration.

\begin{ex}[ (Geometric distribution, 2)]
\label{ex:geom2}
Let the random walk on $\N$ be defined by the transition probabilities
$p_{n,0}=\alpha$ and $p_{n,n+1}=1-\alpha$ for some $0<\alpha<1$. The geometric
distribution $\alpha(1-\alpha)^n$ is invariant (but not reversible) for this random
walk. The Ricci curvature of this random walk is $\alpha$.
\end{ex}

\begin{ex}[ ($\delta$-hyperbolic groups)]
Let $X$ be the Cayley graph of a non-elementary $\delta$-hyperbolic group with respect
to some finite generating set. Let $k$ be a large enough integer
(depending on the group) and consider the random walk
consisting in performing $k$ steps of the simple random walk. Let $x,y\in
X$ with $d(x,y)>2k$. Then $\kappa(x,y)= -2k/d(x,y)+O(1/d(x,y))$.
\end{ex}

Note that $-2k/d(x,y)$ is the smallest possible value for $\kappa(x,y)$,
knowing that the steps of the random walk are bounded by $k$.

\begin{dem}
For $z$ in the ball of radius $k$ around $x$, and $z'$ in the ball of
radius $k$ around $y$, elementary $\delta$-hyperbolic geometry yields $d(z,z')=d(x,y)+d(x,z)+d(y,z')-(y,z)_x-(x,z')_y$ up to
some multiple of $\delta$, where $(\cdot,\cdot)$ denotes the Gromov
product with respect to some basepoint~\cite{GH90}. Since this decomposes
as the sum of a term depending on $z$ only and a term depending on $z'$
only, to compute the transportation distance it is enough to study the expectation of $(y,z)_x$ for $z$ in the
ball around $x$, and likewise for $(x,z')_y$. Knowing that balls have
exponential growth, it is not difficult to see that the expectation of
$(y,z)_x$ is bounded by a constant, whatever $k$, hence the conclusion.

The same argument applies to trees or discrete $\delta$-hyperbolic spaces with a
uniform lower bound on the exponential growth rate of balls.
\end{dem}

\begin{ex}[ (Kac's random walk on orthogonal matrices, after~\cite{Oli})]
Consider the following random walk on the set of $N\times N$ orthogonal
matrices: at each step, a pair of indices $1\leq i<j\leq N$ is selected
at random, an angle $\theta \in [0;2\pi)$ is picked at random, and a
rotation of angle $\theta$ is performed in the coordinate plane $i,j$.
Equip the set of orthogonal matrices
with the Riemannian metric on $\mathrm{SO}(N)$ induced by the
Hilbert--Schmidt inner product $\Tr (a^\ast b)$ on its tangent space.
It is proven in a preprint by Oliveira~\cite{Oli} that this random walk has
curvature $1-\sqrt{1-2/N(N-1)}\sim 1/N^2$.
\end{ex}

This is consistent with the fact that $\mathrm{SO}(N)$ has, as a
Riemannian manifold, a positive Ricci curvature in the usual sense.
However, from the computational point of view, Kac's random walk above is
much nicer than either the Brownian motion or the $\eps$-scale random
walk of Example~\ref{ex:riem}. Oliveira uses this result to prove a new
estimate $O(N^2\ln N)$ for the mixing time of this random walk, nicely improving
on previous estimates $O(N^4 \ln N)$ by Diaconis--Saloff-Coste and
$O(N^{2.5}\ln N)$ by Pak--Sidenko (an easy lower bound is $\Omega(N^2)$),
see~\cite{Oli}.

\begin{ex}[ (Glauber dynamics for the Ising model)]
\label{ex:glauber}
Let $G$ be a finite graph. Consider the
configuration space is $X:=\{-1,1\}^G$ together with the energy function
$U(S):=-\sum_{x\sim y \in G} S(x)S(y)-H\sum_x S(x)$ for $S\in X$, where
$H\in\R$ is the external magnetic field. For some $\beta\geq 0$, equip $X$
with the Gibbs distribution $\mu:=\e^{-\beta U}\!/Z$ where as usual
$Z:=\sum_S \e^{-\beta U(S)}$. The distance between two states is defined
as the number of vertices of $G$ at which their value differ.

For $S\in X$ and $x\in G$, denote by $S_{x+}$ and $S_{x-}$ the states
obtained from $S$ by setting $S_{x+}(x)=+1$ and $S_{x-}(x)=-1$,
respectively. Consider the following random walk on $X$ (known as the
\emph{Glauber dynamics}): at each step, a vertex $x\in G$ is chosen at
random, and a new value for $S(x)$ is picked according to local
equilibrium, i.e.\ $S(x)$ is set to $1$ or $-1$ with probabilities
proportional to $\e^{-\beta U(S_{x+})}$ and $\e^{-\beta U(S_{x-})}$
respectively (note that only the neighbors of $x$ influence the ratio of
these probabilities). The Gibbs distribution is reversible for this
Markov chain.

Then the Ricci curvature of this Markov chain is at least
\[
\frac{1}{\abs{G}}\left(1-v_\text{max}\,\frac{\e^\beta-\e^{-\beta}}{\e^{\beta}+\e^{-\beta}}\right)
\]
where $v_\text{max}$ is the maximal valency of a vertex of $G$. In
particular, if
\[
\beta< \frac12\, \ln \left(\frac{v_\text{max}+1}{v_\text{max}-1}\right)
\]
then curvature is positive. Consequently, the critical $\beta$ is at
least this quantity.
\end{ex}

This estimate for the critical temperature coincides exactly with the one
derived in~\cite{Gri67}; actually our argument generalizes to
non-constant values of the coupling $J_{xy}$ between spins, and the
positive curvature condition exactly amounts to $G(\beta)<1$ in that
paper's notation (\cite{Gri67}, Eq.~(19)), or, equivalently, to
Dobrushin's criterion using a single site.  For comparison, the exact
value of the critical $\beta$ for the Ising model on the regular infinite
tree of valency $v$ is $\frac12 \ln\left(\frac{v}{v-2}\right)$, which
shows asymptotic optimality.

As shown in the rest of this paper, positive curvature implies several
properties, especially, exponential convergence to the equilibrium,
concentration inequalities and a modified logarithmic Sobolev inequality.
I do not know how these results compare to the literature.

Since the argument presented below does not rely on exact solutions but
on quantitative estimates, it is obviously not specific to the Ising
model: the only property we used is that the influence of a vertex on the
local equilibrium of its neighbors is bounded.

\begin{dem}
Using Proposition~\ref{prop:epsgeod}, it is enough to bound Ricci
curvature for pairs states at distance $1$.
Let $S$, $S'$ be two states differing only at $x\in G$. We can
suppose that $S(x)=-1$ and $S'(x)=1$. Let $m_S$ and $m_{S'}$ be the law
of the step of the random walk issuing from $S$ and $S'$ respectively.
We have to prove that the transportation distance between $m_S$ and
$m_{S'}$ is at most
$1-\frac{1}{\abs{G}}\left(1-v_\text{max}\,\frac{\e^{\beta}-\e^{-\beta}}{\e^{\beta}+\e^{-\beta}}\right)$.

The measure $m_S$ decomposes as $m_S=\frac1{\abs{G}} \sum_{y\in G}
m_S^y$, according to the vertex $y\in G$ which is modified by the random
walk, and likewise for $m_{S'}$. To evaluate the transportation distance,
we will compare $m_S^y$ to $m_{S'}^y$.

If the step of the random walk consists in modifying the value of $S$ at
$x$ (which occurs with probability $1/\abs{G}$), then the resulting state
has the same law for $S$ and $S'$, i.e.\ $m_S^x=m_{S'}^x$. Thus in this
case the transportation distance is $0$ and the contribution to Ricci
curvature is $1\times \frac1{\abs{G}}$.

If the step consists in modifying the value of $S$ at some point $y$ in
$G$ not adjacent to $x$, then the value at $x$ does not influence local
equilibrium at $y$, and so $m_S^y$ and $m_{S'}^y$ are identical except at
$x$.  So in this case the distance is $1$ and the contribution to Ricci
curvature is $0$.

Now if the step consists in modifying the value of $S$ at some point
$y\in G$ adjacent to $x$ (which occurs with probability $v_x/\abs{G}$
where $v_x$ is the valency of $x$), then the value at $x$ does influence
the law of the new value at $y$, by some amount which we now evaluate.
The final distance between the two laws will be this amount plus $1$ ($1$
accounts for the difference at $x$), and the contribution to Ricci
curvature will be negative.

Let us now evaluate this amount more precisely. Let $y\in G$ be adjacent
to $x$. Set $a=\e^{-\beta U(S_{y+})}/\e^{-\beta U(S_{y-})}$. The step of
the random walk consists in setting $S(y)$ to $1$ with probability
$\frac{a}{a+1}$, and to $-1$ with probability $\frac1{a+1}$. Setting likewise
$a'=\e^{-\beta U(S'_{y+})}/\e^{-\beta U(S'_{y-})}$ for $S'$, we are left to
evaluate the distance between the distributions on $\{-1,1\}$ given by
$\left(\frac{a}{a+1};\frac1{a+1}\right)$ and
$\left(\frac{a'}{a'+1};\frac1{a'+1}\right)$. It is immediate to check,
using the definition of the energy $U$, that $a'=\e^{4\beta} a$. Then, a
simple computation shows that the distance between these two
distributions is at most
$\frac{\e^{\beta}-\e^{-\beta}}{\e^{\beta}+\e^{-\beta}}$. This value
is actually achieved when $y$ has odd valency, $H=0$ and switching the value at
$x$ changes the majority around $y$. (Our argument is suboptimal here
when valency is even---a more precise estimation yields the absence of a
phase transition in dimension $1$.)

Combining these different cases yields the desired curvature evaluation.
To convert this into an evaluation of the critical $\beta$, reason as
follows: Magnetization, defined as $\frac1{\abs{G}}\sum_{x\in G} S(x)$,
is a $\frac{1}{\abs{G}}$-Lipschitz function of the state. Now let $\mu_0$ be
the Gibbs measure without magnetic field, and $\mu_h$ the Gibbs measure
with external magnetic field $h$.
Use the Glauber dynamics with
magnetic field $h$, but starting with an initial state picked under
$\mu_0$; Cor.~\ref{cor:avlip} yields that the magnetization under $\mu_h$ is
controlled by $\frac1{\abs{G}}\,\Td(\mu_0,\mu_0\ast m)/\kappa$ where
$\kappa$ is the Ricci curvature, and $\Td(\mu_0,\mu_0\ast m)$ is the
transportation distance between the Gibbs measure $\mu_0$ and the measure
obtained from it after one step of the Glauber dynamics with magnetic
field $h$; reasoning as above this transportation distance is easily bounded by
$\frac1{\abs{G}}\,\frac{\e^{\beta h}-\e^{-\beta h}}{\e^{\beta h}+\e^{-\beta
h}}$, so that the derivative of the magnetization w.r.t.\ $h$ stays
bounded
when $\abs{G}\to\infty$. (Compare Eq.~(22) in~\cite{Gri67}.)
\end{dem}

\medskip

More examples can be found in Sections~\ref{sec:poisson} (binomial and
Poisson distributions), \ref{sec:continuoustime} ($M/M/\infty$ queues
and generalizations) and~\ref{sec:exp} (geometric distributions on $\N$,
exponential distributions on $\R^N$).

\subsection{Overview of the results}

\paragraph{Notation for random walks.} Before presenting the main
results, we need some more quantites related to the local behavior of the
random walk: the \emph{jump}, which will help control the diameter of the
space, and the \emph{spread}, which is the analogue of a diffusion
constant and will help control concentration properties. Moreover, we
define a notion of local dimension. The larger the dimension, the
better for concentration of measure.
 

\begin{defi}[ (Jump, spread, dimension)]
\label{def:jsd}
Let the \emph{jump} of the random walk at $x$ be
\[
J(x):=\E_{m_x} d(x,\cdot)=\Td(\delta_x,m_x)
\]

Let the \emph{spread} of the random walk at $x$ be
\[
\sigma(x):=\left(\frac12 \iint d(y,z)^2 \,\d m_x(y) \,\d m_x(z)\right)^{1/2} 
\]
and, if $\nu$ is a invariant distribution, let
\[
\sigma:=\norm{\sigma(x)}_{L^2(X,\nu)}
\]
be the average spread.

Let also
$
\sigma_\infty(x):=\frac12 \diam \Supp m_x
$
and
$
\sigma_\infty:=\sup \sigma_\infty(x)
$.

Let the \emph{local dimension} at $x$ be
\[
n_x:=\frac{\sigma(x)^2}{\sup \{\Var_{m_x} f, f\,
1\text{-Lipschitz}\}}
\]
and finally $n:=\inf_x n_x$.
\end{defi}

\paragraph{About this definition of dimension.}
Obviously $n_x\geq 1$. For the discrete-time Brownian motion on a
$N$-dimensional Riemannian manifold, one has $n_x\approx N$ (see the end of
Section~\ref{sec:riemproof}). For the simple
random walk on a graph, $n_x\approx 1$. This definition of
dimension amounts to saying that in a space of dimension $n$, the typical
variations of a ($1$-dimensional) Lipschitz function are $1/\sqrt{n}$
times the typical distance between two points.  This is the case in the
sphere $S^n$, in the Gaussian measure on $\R^n$, and in the discrete cube
$\{0,1\}^n$. So generally one could define the ``statistical dimension'' of
a metric measure space $(X,d,\mu)$ by this formula i.e.
\[
\mathrm{StatDim}(X,d,\mu):=\frac{\frac12 \iint d(x,y)^2
\,\d\mu(x)\d\mu(y)}{\sup \{\Var_{\mu} f, f\,
1\text{-Lipschitz}\}}
\]
so that for each $x\in X$ the local dimension of $X$ at $x$ is 
$n_x=\mathrm{StatDim}(X,d,m_x)$.
With this definition, $\R^N$ equipped with a Gaussian measure has statistical
dimension $N$ and local dimension $\approx N$, whereas the discrete cube
$\{0,1\}^N$ has statistical dimension $\approx N$ and local dimension $\approx 1$.

\bigskip

We now turn to the description of the main results of the paper.

\paragraph{Elementary properties.} In Section~\ref{sec:basic} are
gathered some straightforward results.

First, we prove (Proposition~\ref{prop:epsgeod}) that in an
$\eps$-geodesic space, it is enough to get a lower bound on $\kappa(x,y)$
for points $x,y$ with $d(x,y)\leq \eps$, to get a lower bound on $\kappa$
for all pairs of points. This is simple yet very useful: indeed in the various
graphs given above as examples, it was enough to compute the Ricci
curvature for neighbors.

Second, we prove equivalent characterizations of having Ricci curvature
uniformly bounded from below: A space satisfies $\kappa(x,y)\geq \kappa$
if and only if the random walk operator is $(1-\kappa)$-contracting on
the space of probability measures equipped with the transportation
distance (Proposition~\ref{prop:contracting}), and if and only if the
random walk operator acting on Lipschitz functions contracts the
Lipschitz norm by $(1-\kappa)$ (Proposition~\ref{prop:lipcontracting}).
An immediate corollary of the contracting property for probability
measures is the existence of a unique invariant distribution when
$\kappa>0$.

The property of contraction of the Lipschitz norm implies, in the
reversible case, that the spectral gap of the Laplacian operator
associated with the random walk is at least $\kappa$; this can be seen as
a generalization of Lichnerowicz' theorem, and provides sharp estimates
of the spectral gap in several examples.

In analogy with the Bonnet--Myers theorem, we prove that if Ricci
curvature is bounded below by $\kappa>0$, then the diameter of the space
is at most $2\sup_x J(x)/\kappa$ (Proposition~\ref{prop:weakBM}). In case
$J$ is unbounded, we can evaluate instead the average distance to a given
point $x_0$ under the invariant distribution $\nu$
(Proposition~\ref{prop:meandiam}); namely, $\int d(x_0,y)\,\d\nu(y)\leq
J(x_0)/\kappa$. In particular we have $\int d(x,y)\,\d\nu(x)\d\nu(y)\leq
2\inf J/\kappa$. These are $L^1$ versions of the Bonnet--Myers theorem
rather than generalizations: from the case of manifolds one would expect
$1/\sqrt{\kappa}$ instead of $1/\kappa$. Actually this $L^1$ version is
sharp in all our examples except Riemannian manifolds; in
Section~\ref{sec:strongBM} we investigate additional conditions for an $L^2$
version of the Bonnet--Myers theorem to hold.

Let us also mention two elementary constructions preserving positive
curvature, namely, superposition and $L^1$ tensorization
(Propositions~\ref{prop:superposition} and~\ref{prop:tensorization}).

\paragraph{Concentration results.} Basically, if Ricci curvature is
bounded below by $\kappa>0$, then the invariant distribution satisfies
concentration results with variance $\sigma^2/n\kappa$ (up to some
constant factor). This estimate is often sharp, as discussed in
Section~\ref{sec:examplesrevisited} where we revisit some of the
examples.

However, the type of concentration (Gaussian, exponential, or $1/t^2$)
depends on further local assumptions: indeed, just as in the central
limit theorem, positive Ricci curvature can only carry at the global
scale what is already true at the local scale.  Without further
assumptions, one only gets that the maximal variance of a $1$-Lipschitz
function is at most $\sigma^2/n\kappa$, hence concentration like
$\sigma^2/n\kappa t^2$ (Proposition~\ref{prop:varlip}). If we make the
further assumption that the support of the measures $m_x$ is uniformly
bounded (i.e.\ $\sigma_\infty<\infty$), then we get mixed
Gaussian-then-exponential concentration, with variance $\sigma^2/n\kappa$
(Theorem~\ref{thm:gaussconc}). The width of the Gaussian window depends
on $\sigma_\infty$, and on the rate of variation of the spread
$\sigma(x)^2$.

For the case of Riemannian manifolds, simply taking smaller and smaller
steps for the random walks makes the width of the Gaussian window tend to
infinity, so that we recover Gaussian concentration as in the
Lévy--Gromov or Bakry--Émery context. However, for lots of discrete
examples, the Gaussian-then-exponential behavior is genuine. Examples
where tails are Poisson-like (binomial distribution, $M/M/\infty$ queues)
or exponential are given in Sections~\ref{sec:poisson}
to~\ref{sec:unstability}.

We also get concentration results for the finite-time distributions
$m_x^{\ast k}$ (Remark~\ref{rem:finitetime}).

\paragraph{Log-Sobolev inequality.} Using a suitable non-local notion of
norm of the gradient, we are able to mimic the proof by Bakry and Émery
of a logarithmic Sobolev inequality for the invariant distribution.  The
gradient we use (Definition~\ref{def:gradient}) is
$(Df)(x):=\sup_{y,z}\frac{\abs{f(y)-f(z)}}{d(y,z)}\,\exp(-\lambda
d(x,y)-\lambda d(x,z))$. This is a kind of ``semi-local'' Lipschitz
constant for $f$. Typically the value of $\lambda$ can be taken large at
the ``macroscopic'' level; for Riemannian manifolds, taking smaller and
smaller steps for the random walk allows to take $\lambda\to \infty$ so
that we recover the usual gradient for smooth functions.

The inequality takes the form $\Ent f\leq C \int (Df)^2/f \,\d\nu$
(Theorem~\ref{thm:logsob}).  The main tool of the proof is the
contraction relation $D(M\!f)\leq (1-\kappa/2) M(Df)$ where $M$ is the
random walk operator (Proposition~\ref{prop:localcontrol}).

That the gradient is non-local, with a maximal possible value of
$\lambda$, is consistent with the possible occurrence of non-Gaussian
tails.

\paragraph{Exponential concentration and non-negative curvature.} The
simplest example of a Markov chain with zero Ricci curvature is the
simple random walk on $\N$ or $\Z$, for which there is no invariant
distribution. However, we show that if furthermore there is a ``locally
attracting'' point, then non-negative Ricci curvature implies exponential
concentration. The main examples are the geometric distribution on $\N$,
and the exponential distribution $\e^{-\abs{x}}$ on $\R^n$ associated with
the stochastic differential equation $\d X_t=\d
B_t-\frac{X_t}{\abs{X_t}}\,\d t$. In both cases we recover correct orders
of magnitude.

\paragraph{Gromov--Hausdorff topology.} One advantage of our definition
is that it involves only combinations of the distance function, and no
derivatives, so that it is more or less impervious to deformations of the
space. In Section~\ref{sec:GH} we show that Ricci curvature is continuous
for Gromov--Hausdorff convergence of metric spaces (suitably reinforced,
of course, so that the random walk converges as well), so that having
non-negative curvature is a closed property. We also suggest a loosened
definition of Ricci curvature, requiring that $\Td(m_x,m_y)\leq
(1-\kappa)d(x,y)+\delta$ instead of $\Td(m_x,m_y)\leq (1-\kappa)d(x,y)$.
With this definition, positive curvature becomes an \emph{open} property,
so that a space close to one with positive curvature has positive
curvature. Properties of this loose version will be investigated in
another paper.

\section{Elementary properties}

\label{sec:basic}

\subsection{Geodesic spaces}

The idea behind curvature is to use local properties to derive global
ones.  We give here a simple proposition expressing that in near-geodesic
spaces, such as graphs or manifolds, it is enough to check positivity of
Ricci curvature for nearby points.

\begin{prop}
\label{prop:epsgeod}
Suppose that $(X,d)$ is $\eps$-geodesic in the sense that for any two
points $x,y \in X$, there exists an integer $n$ and a sequence
$x_0=x,x_1,\ldots,x_n=y$ such that $d(x_i,x_{i+1})\leq \eps$ and
$d(x,y)=\sum d(x_i,x_{i+1})$.

Then, if $\kappa(x,y)\geq \kappa$ for any pair of points with $d(x,y)\leq
\eps$, then $\kappa(x,y)\geq \kappa$ for any pair of points $x,y\in X$.
\end{prop}

\begin{dem}
Since $\Td$ is a distance, one has $\Td(m_x,m_y)\leq
\sum \Td(m_{x_i},m_{x_{i+1}})\leq
(1-\kappa) \sum d(x_i,x_{i+1})$.
\end{dem}

\subsection{Contraction on the space of probability measures}

Let $\Pm(X)$ by the space of all probability measures $\mu$ on $X$ with
finite first moment, i.e.\ for some (hence any) $o\in X$, $\int
d(o,x)\,\d\mu(x)<\infty$. On $\Pm(X)$, the transportation distance $\Td$ is
finite, so that it is actually a distance.

Let $\mu$ be a probability measure on $X$ and define the measure
\[
\mu\ast m:=\int_{x\in X} \d\mu(x)\, m_x 
\]
which is the image of $\mu$ by the random walk. (It may or
may not belong to $\Pm(X)$.)

The following proposition also appears in~\cite{DGW04} (in the proof of
Proposition~2.10) and in~\cite{Oli}.

\begin{prop}
\label{prop:contracting}
Let $(X,d,m)$ be a metric space with a random walk. Let $\kappa\in \R$. Then the we have
$\kappa(x,y)\geq \kappa$ for all $x,y\in X$, if and only if for any two
probability distributions $\mu,\mu'\in \Pm(X)$ one has
\[
\Td(\mu\ast m,\mu'\ast m)\leq (1-\kappa) \Td(\mu,\mu')
\]

Moreover in this case, if $\mu\in \Pm(X)$ then $\mu\ast m\in \Pm(X)$.
\end{prop}

\begin{dem}
First, suppose that convolution with $m$ is contracting in $\Td$
distance. For some $x,y\in X$, let $\mu=\delta_x$ and $\mu'=\delta_y$ be
the Dirac measures at $x$ and $y$. Then by definition $\delta_x\ast
m=m_x$ and likewise for $y$, so that $\Td(m_x,m_y)\leq
(1-\kappa)\Td(\delta_x,\delta_y)=(1-\kappa)d(x,y)$ as required.

The converse is more difficult to write than to understand.
For each pair $(x,y)$ let $\xi_{xy}$ be a coupling (i.e.\ a measure on $X\times
X$) between $m_x$ and $m_y$ witnessing for $\kappa(x,y)\geq \kappa$.
According to Corollary~5.22 in~\cite{Vil}, we can choose $\xi_{xy}$ to
depend measurably on the pair $(x,y)$.
Let $\Xi$ be a coupling between $\mu$ and $\mu'$ witnessing for
$\Td(\mu,\mu')$. Then $\int_{X\times X} \d\Xi(x,y)\,\xi_{xy}$ is a coupling
between $\mu\ast m$ and $\mu'\ast m$ and so
\begin{eqnarray*}
\Td(\mu\ast m,\mu'\ast m) &\leq&
\int_{x,y}d(x,y) \,\d\!\left\{\int_{x',y'} \d\Xi(x',y')\,
\xi_{x'y'}\right\}\!(x,y)
\\
&=&\int_{x,y,x',y'} \d\Xi(x',y') \, \d\xi_{x'y'}(x,y)\,d(x,y)
\\
&\leq& \int_{x',y'} \d\Xi(x',y') \, d(x',y')(1-\kappa(x',y'))
\\&\leq& (1-\kappa) \Td(\mu,\mu')
\end{eqnarray*}
by the Fubini theorem applied to $d(x,y) \,\d\Xi(x',y')
\,\d\xi_{x',y'}(x,y)$.

To see that in this situation $\Pm(X)$ is preserved by the random walk,
fix some origin $o\in X$ and note that for any $\mu\in \Pm(X)$, the first
moment of $\mu\ast m$ is $\Td(\delta_o,\mu \ast m)\leq
\Td(\delta_o,m_o)+\Td(m_o,\mu_\ast m)\leq
\Td(\delta_o,m_o)+(1-\kappa)\Td(o,\mu)$. Now $\Td(o,\mu)<\infty$ by
assumption, and $\Td(\delta_o,m_o)<\infty$ by our definition of random
walks (Definition~\ref{def:space}).
\end{dem}

As an immediate consequence of this contracting property we get:

\begin{cor}
Suppose that $\kappa(x,y)\geq\kappa>0$ for any two distinct $x,y\in X$.
Then the random walk has a unique invariant distribution $\nu\in \Pm(X)$.

Moreover, for any probability measure $\mu\in\Pm(X)$, the sequence $\mu\ast
m^{\ast n}$ tends exponentially fast to $\nu$ in $\Td$ distance. Namely
\[
\Td(\mu\ast m^{\ast n},\nu)\leq (1-\kappa)^n \Td(\mu,\nu)
\]
and in particular
\[
\Td(m_x^{\ast n},\nu)\leq (1-\kappa)^n J(x)/\kappa
\]
\end{cor}

The last assertion follows by taking $\mu=\delta_x$ and noting that
$J(x)=\Td(\delta_x,m_x)$ so that $\Td(\delta_x,\nu)\leq
\Td(\delta_x,m_x)+\Td(m_x,\nu)\leq J(x)+(1-\kappa)\Td(\delta_x,\nu)$,
hence $\Td(\delta_x,\nu)\leq J(x)/\kappa$.

Another interesting corollary is the following, which allows to estimate
the average of a Lipschitz function under the invariant measure, knowing
some of its values. This is useful in concentration theorems, to get
bounds not only on the deviations from the average, but on what the
average actually is.

\begin{cor}
\label{cor:avlip}
Suppose that $\kappa(x,y)\geq\kappa>0$ for any two distinct $x,y\in X$.
Let $\nu$ be the invariant distribution.

Let $f$ be a $1$-Lipschitz function. Then, for any distribution $\mu$,
one has $\abs{\E_\nu f-\E_\mu f}\leq \Td(\mu,\mu\ast m)/\kappa$.

In particular, for any $x\in X$ one has $\abs{f(x)-\E_\nu f}\leq
J(x)/\kappa$.
\end{cor}

\begin{dem}
One has $\Td(\mu\ast m,\nu)\leq (1-\kappa)\Td (\mu,\nu)$. Since by the
triangle inequality, $\Td(\mu\ast m, \nu)\geq
\Td(\mu,\nu)-\Td(\mu,\mu\ast m)$, one gets $\Td(\mu,\nu)\leq
\Td(\mu,\mu\ast m)/\kappa$. Now if $f$ is a $1$-Lipschitz function, for
any two distributions $\mu$,$\mu'$ one has $\abs{\E_\mu f-\E_{\mu'}f}\leq
\Td(\mu,\mu')$ hence the result.

The last assertion is simply the case when $\mu$ is the Dirac measure at
$x$.
\end{dem}

\subsection{$L^1$ Bonnet--Myers theorems}
\label{sec:weakBM}

We now give a weak analogue of the Bonnet--Myers theorem.  This result
shows in particular that positivity of Ricci curvature is a much stronger
property than some spectral gap bound: there is no Ricci curvature
analogue of a family of expanders.

\begin{prop}[ ($L^1$ Bonnet--Myers)]
\label{prop:weakBM}
Suppose that $\kappa(x,y)\geq \kappa>0$ for all
$x,y\in X$. Then for any $x,y\in X$ one has
\[
d(x,y)\leq \frac{J(x)+J(y)}{\kappa(x,y)}
\]
and in particular
\[
\diam X \leq \frac{2\sup_x J(x)}{\kappa}
\]
\end{prop}

\begin{dem}
Let $d=d(x,y)$. By assumption we have $\Td(m_x,m_y)\leq d(1-\kappa)$. By
definition we have $\Td(m_x,\delta_x)=J(x)$ and
$\Td(m_y,\delta_y)=J(y)$. So $d\leq J(x)+J(y)+d(1-\kappa)$.
\end{dem}

This result is not sharp at all for Brownian motion in Riemannian
manifolds (since $J\approx \eps$ and $\kappa\approx \eps^2 \Ric/N$, it
fails by a factor $1/\eps$ compared to the Bonnet--Myers
theorem!), but is sharp in many other examples.

For the discrete cube $X=\{0,1\}^N$ (Example~\ref{ex:cube} above), one has
$J=1/2$ and $\kappa=1/N$, so we get $\diam X\leq N$ which is the exact
value.

For the discrete Ornstein--Uhlenbeck process (Example~\ref{ex:discrOU}
above) one has $J=1/2$ and $\kappa=1/2N$, so we get $\diam X\leq 2N$
which once more is the exact value.

For the continuous Ornstein--Uhlenbeck process on $\R$
(Example~\ref{ex:OU} with $N=1$), the diameter is infinite, consistently
with the fact that $J$ is unbounded. If we restrict the process to some
large interval $[-R;R]$ with $R\gg s/\sqrt{\alpha}$ (e.g.\ by reflecting
the Brownian part), then $\sup J\sim \alpha R\deltat$ on this interval,
and $\kappa = (1-\e^{\alpha\deltat})\sim \alpha\deltat$ so that the
diameter is bounded by $2R$, which is correct.

These examples show that one cannot replace $J/\kappa$ with
$J/\sqrt{\kappa}$ in this result (as could be expected from the example
of Riemannian manifolds).  In fact, Riemannian manifolds seem to be the
only simple example where there is a diameter bound behaving like
$1/\sqrt{\kappa}$. In Section~\ref{sec:strongBM} we investigate
conditions under which an $L^2$ version of the Bonnet--Myers theorem
holds.

\bigskip

In case $J$ is not bounded, we can estimate instead the ``average''
diameter $\int d(x,y)\, \d\nu(x) \d\nu(y)$ under the invariant
distribution $\nu$.  This estimate will prove very useful in several
examples, to get bounds on the average of $\sigma(x)$ in cases where
$\sigma(x)$ is unbounded but controlled by the distance to some
``origin'' (see e.g.\ Sections~\ref{sec:continuoustime}
and~\ref{sec:unstability}).

\begin{prop}[ (Average $L^1$ Bonnet--Myers)]
\label{prop:meandiam}
Suppose that $\kappa(x,y)\geq\kappa>0$ for any two distinct $x,y\in X$.
Then for any $x\in X$,
\[
\int_X d(x,y) \,\d\nu(y)\leq \frac{J(x)}{\kappa}
\]
and so
\[
\int_{X\times X} d(x,y) \,\d\nu(x)\,\d\nu(y)\leq \frac{2 \inf_x J(x)}{\kappa}
\]
\end{prop}

\begin{dem}
%
The first assertion follows from Corollary~\ref{cor:avlip} with
$f=d(x,\cdot)$.

For the second assertion, choose an $x_0$ such that $J(x_0)$ is arbitrarily close to $\inf
J$, and write
\begin{eqnarray*}
\int_{X\times X} d(y,z) \,\d\nu(y) \,\d\nu(z)
&\leq &
\int_{X\times X} (d(y,x_0)+d(x_0,z)) \,\d\nu(y) \,\d\nu(z)
\\&=& 2\Td (\delta_{x_0},\nu)\leq 2J(x_0)/\kappa
\end{eqnarray*}
which ends the proof.
\end{dem}

\subsection{Two constructions}

Here we describe two very simple constructions which trivially preserve
positive curvature, namely, superposition and $L^1$ tensorization.

Superposition states that if we are given two random walks on the same
space and construct a new one by, at each step, tossing a coin and
deciding to follow either one random walk or the other, then the Ricci
curvatures mix nicely.

\begin{prop}[ (Superposition)]
\label{prop:superposition}
Let $X$ be a metric space equipped with a family $(m^{(i)})$ of random
walks.
Suppose that for each $i$, the Ricci curvature of $m^{(i)}$ is at least
$\kappa_i$. Let $(\alpha_i)$ be a family of non-negative real numbers
such that $\sum \alpha_i=1$.
Define a random walk $m$ on $X$ by $m_x:=\sum \alpha_i m^{(i)}_x$. Then the
Ricci curvature of $m$ is at least $\sum\alpha_i\kappa_i$.
\end{prop}

\begin{dem}
Let $x,y\in X$ and for each $i$ let $\xi_i$ be a couplings between
$m^{(i)}_x$ and $m{(i)}_y$.
Then $\sum\alpha_i\xi_i$ is a
coupling between $\sum \alpha_i m^{(i)}_x$ and $\sum\alpha_i m^{(i)}_y$, so that 
\begin{align*}
\Td(m_x,m_y)&\leq \sum \alpha_i \,\Td\left(m^{(i)}_x,m^{(i)}_y\right)
\\& \leq
\sum \alpha_i(1-\kappa_i) d(x,y)
\\&
=\left(1-\sum\alpha_i\kappa_i\right)d(x,y)
\end{align*}

Note that the coupling above, which consists in sending each $m^{(i)}_x$
to $m^{(i)}_y$, has no reason to be optimal, so that in general equality
does not hold.
\end{dem}

Tensorization states that if we perform a random walk in a product space
by deciding at random, at each step, to move in one or the other
component, then positive curvature is preserved.

\begin{prop}[ ($L^1$ tensorization)]
\label{prop:tensorization}
Let $(X_1,\ldots,X_k)$ be a finite family of metric spaces equipped with a family of random
walks $(m^{(1)},\ldots,m^{(k)})$.
Let $X$ be the product of the spaces $X_i$,
equipped with the distance
$\sum d_i$.
Let $(\alpha_i)$ be a family of non-negative real numbers such that $\sum
\alpha_i=1$.
Consider the random walk on $X$ defined by
\[
m_{(x_i)}:=\sum \alpha_i\,\, \delta_{x_1}\otimes\cdots\otimes
m_{x_i}\otimes \cdots \otimes \delta_{x_k}
\]

Suppose that for each $i$, the Ricci curvature of $m^{(i)}$ is at least
$\kappa_i$.
Then the Ricci curvature of $m$ is at least
$\inf \alpha_i \kappa_i$.
\end{prop}

For example, this allows for a very short proof that the curvature of the
lazy random walk on the discrete cube $\{0,1\}^N$ is $1/N$
(Example~\ref{ex:cube}). Indeed, it is the $N$-fold product of the random
walk on $\{0,1\}$ which sends each point to the equilibrium distribution
$(1/2,1/2)$, hence is of curvature $1$.

The case when some $\alpha_i$ is equal to $0$ shows why the Ricci
curvature is given by an infimum: indeed, if $\alpha_i=0$ then the
corresponding component never gets mixed, hence curvature cannot be
positive (unless this component is reduced to a single point).

Here the statement is restricted to a finite product for the following
technical reasons: First, to define the $L^1$ product of an infinite
family, a basepoint has to be chosen. Second, in order for the formula
above to define a random walk with finite first moment (see
Definition~\ref{def:space}), some uniform assumption on the first moments
of the $m^{(i)}$ is needed.

\begin{dem}
For $x\in X$ let $\tilde{m}^{(i)}_x$ stand for
$\delta_{x_1}\otimes\cdots\otimes
m_{x_i}\otimes \cdots \otimes \delta_{x_k}$.

Let $x=(x_i)$ and $y=(y_i)$ be two points in $X$. 
Then
\begin{align*}
\Td(m_x,m_y)&\leq \sum
\alpha_i\,\Td\left(\tilde{m}^{(i)}_x,\tilde{m}^{(i)}_y\right)
\\&\leq
\sum \alpha_i\left(\Td\left(m^{(i)}_x,m^{(i)}_y\right)+\sum_{j\neq i}
d_j(x_j,y_j)\right)
\\&\leq
\sum \alpha_i\left((1-\kappa_i)d_i(x_i,y_i)+\sum_{j\neq i}
d_j(x_j,y_j)\right)
\\&=
\sum \alpha_i\left(-\kappa_id_i(x_i,y_i)+\sum
d_j(x_j,y_j)\right)
\\&=
\sum d_i(x_i,y_i)-\sum\alpha_i\kappa_id_i(x_i,y_i)
\\&\leq (1-\inf \alpha_i\kappa_i)\,d(x,y)
\end{align*}
\end{dem}

\subsection{Lipschitz functions and spectral gap}

\begin{defi}[ (Averaging operator, Laplacian)]
For $f\in L^2(X,\nu)$ let the averaging operator $M$ be
\[
M\!f(x):=\int_y f(y) \,\d m_x(y)
\]
and let $\Delta:=M-\Id$.
\end{defi}

(This is the layman's convention for the sign of the Laplacian, i.e.\
$\Delta=\frac{d^2}{dx^2}$ on $\R$, so that on a Riemannian manifold $\Delta$
is a negative operator.)

The following proposition also appears in~\cite{DGW04} (in the proof of
Proposition~2.10).

\begin{prop}
\label{prop:lipcontracting}
Let $(X,d,m)$ be a random walk on a metric space. Let $\kappa\in \R$.

Then the Ricci curvature of $X$ is at least $\kappa$, if and only if, for
every $k$-Lipschitz function
$f:X\to \R$, the function $M\!f$ is $k(1-\kappa)$-Lipschitz.
\end{prop}

\begin{dem}
First, suppose that the Ricci curvature of $X$ is at least $\kappa$.
Then
we have
\begin{align*}
M\!f(y)-M\!f(x)&=\int_z f(y+z)-f(x+z)
\\&\leq k \int_z d(x+z,y+z)
\\&= k d(x,y) (1-\kappa(x,y))
\end{align*}

Conversely, suppose that whenever $f$ is $1$-Lipschitz, $M\!f$ is
$(1-\kappa)$-Lipschitz. The duality theorem for transportation distance
(Theorem~1.14 in~\cite{Vil03}) states that
\begin{align*}
\Td(m_x,m_y)&=\sup_{f\text{ $1$-Lipschitz}} \int f \,\d(m_x-m_y)
\\&=\sup_{f\text{ $1$-Lipschitz}} M\!f(x)-M\!f(y)
\\&\leq (1-\kappa) d(x,y)
\end{align*}
\end{dem}

Let $\nu$ be an invariant distribution of the random walk. 
Consider the space $L^2(X,\nu)/\{\text{const}\}$ equipped with the norm
$\norm{f}^2_{L^2(X,\nu)/\{\text{const}\}}:=\norm{f-\int f
\d\nu}^2_{L^2(X,\nu)}$ so that
\[
\norm{f}^2_{L^2(X,\nu)/\{\text{const}\}}=
\Var_\nu f = \frac12
\int_{X\times X} (f(x)-f(y))^2 \,\d\nu(x) \,\d\nu(y)
\]
The operators $M$
and $\Delta$ are self-adjoint in $L^2(X,\nu)$ if and only if $\nu$ is
reversible for the random walk.

It is easy to check, using associativity of variances, that
\[
\Var_\nu f=\int \Var_{m_x} f \,\d\nu(x)+\Var_\nu M\!f
\]
so that $\norm{M\!f}_2\leq \norm{f}_2$. It is also clear that
$\norm{M\!f}_\infty\leq \norm{f}_\infty$.

Usually, spectral gap properties for $\Delta$ are expressed in the space
$L^2$. The proposition above only implies that the spectral radius of the
operator $M$ acting on $\text{Lip}(X)/\{{\text{const}\}}$ is at most
$(1-\kappa)$.  In general it is not true that a bound for the spectral
radius of an operator on a dense subspace of a Hilbert space implies a
bound for the spectral radius on the whole space. This holds, however,
when the operator is self-adjoint or when the Hilbert space is
finite-dimensional.

\begin{prop}
Let $(X,d,m)$ be metric space with random walk, with
invariant distribution $\nu$. Suppose that the Ricci curvature of
$X$ is at least $\kappa>0$ and that $\sigma<\infty$.
Suppose that $\nu$ is reversible, or that $X$ is finite.

Then the spectral radius of the averaging
operator acting on $L^2(X,\nu)/\{\text{const}\}$ is at most $1-\kappa$.
\end{prop}

\begin{dem}
First, if $X$ is finite then Lipschitz functions coincide with $L^2$
functions, so that there is nothing to prove. So we suppose that $\nu$ is
reversible, i.e.\ $M$ is self-adjoint.

Let $f$ be a $k$-Lipschitz function. Proposition~\ref{prop:varlip}
below implies that Lipschitz functions belong to $L^2$ and that the
Lipschitz norm controls the $L^2$ norm. (This is where we use that
$\sigma<\infty$.)

Since $M^t \!f$ is
$k(1-\kappa)^t$-Lipschitz one gets $\Var M^t \!f\leq C k^2 (1-\kappa)^{2t}$
for some constant $C$ so that $\lim_{t\to \infty} (\sqrt{\Var M^t \!f})^{1/t}\leq
(1-\kappa)$. Now Lipschitz functions are dense in $L^2(X,\nu)$. Since $M$
is bounded and self-adjoint, its spectral radius is at most $1-\kappa$.
\end{dem}

\begin{cor}
\label{cor:poincare}
Let $(X,d,m)$ be an ergodic random walk on a metric space, with
invariant distribution $\nu$. Suppose that the Ricci curvature of
$X$ is at least $\kappa>0$ and that $\sigma<\infty$.
Suppose that $\nu$ is reversible.

Then the smallest eigenvalue of $-\Delta$ on
$L^2(X,\nu)/\{\text{const}\}$ is at least $\kappa$.

Moreover the following discrete Poincaré inequalities are satisfied for $f\in
L^2(X,\nu)$:
\[
\Var_\nu f\leq \frac{1}{\kappa(2-\kappa)}\,\int \Var_{m_x}f \,\d\nu(x)
\]
and
\[
\Var_\nu f\leq \frac{1}{2\kappa}\,\iint (f(y)-f(x))^2\,\d\nu(x)
\,\d m_x(y)
\]
\end{cor}

\begin{dem}
These are rewritings of the inequalities $\Var_\nu M\!f\leq (1-\kappa)^2
\Var_\nu f$ and $\langle f,M\!f\rangle_{L^2(X,\nu)/\{\text{const}\}}\leq
(1-\kappa)\Var_\nu f$, respectively.
\end{dem}

The quantities $\Var_{m_x} f$ and $\frac12 \int (f(y)-f(x))^2\,\d
m_x(y)$  are two possible definitions of $\norm{\nabla \!f(x)}^2$ in a
discrete setting. Though the latter is more common, the former is
preferable when the support of $m_x$ can be far away from $x$ and cancels
out the ``drift''. Moreover one always has $\Var_{m_x}f\leq \int
(f(y)-f(x))^2\,\d m_x(y)$, so that the first form is generally sharper
(note that since $\kappa\leq 1$ one has $1/\kappa(2-\kappa)\leq
1/\kappa$).

Reversibility is really needed here to turn an estimate of the spectral
radius of $M$ into an inequality between the norms of $M\!f$ and $f$,
using that $M$ is self-adjoint.  When the random walk is not reversible,
a version of the Poincaré inequality with a non-local gradient still
holds (Theorem~\ref{thm:logsob}).

\bigskip

Let us compare this result to Lichnerowicz' theorem in the case of the
random walk at scale $\eps$ on an $N$-dimensional Riemannian manifold
with positive Ricci curvature.
The operator $\Delta$ associated with the random walk is the difference
between the mean value of a function on a ball of radius $\eps$, and its
value at the center of the ball: when $\eps\to 0$ this behaves like
$\frac{\eps^2}{2(N+2)}$ times the usual Laplacian, by taking the average
on the ball of the Taylor expansion of $f$. Meanwhile, we saw
(Example~\ref{ex:riem}) that $\kappa\sim \frac{\eps^2}{2(N+2)}\,\inf
\Ric$, where $\inf \Ric$ is the largest $K$ such that $\Ric(v,v)\geq K$
for all unit tangent vectors $v$. Note that both scaling factors are the
same. On the other hand the Lichnerowicz theorem states that the smallest
eigenvalue of the usual Laplacian is $\frac{N}{N-1}\inf \Ric$. So we miss
the $\frac{N}{N-1}$ factor, but otherwise get the correct order of
magnitude.

Second, let us test this corollary for the discrete cube of
Example~\ref{ex:cube}. In this case the eigenbase of the discrete
Laplacian is well-known (characters, or Fourier/Walsh transform), and the
spectral radius of the lazy random walk is exactly $1-1/N$. Since the
Ricci curvature $\kappa$ is $1/N$, the value given in the proposition is
sharp.

Third, consider the Ornstein--Uhlenbeck process on $\R$, as in
Example~\ref{ex:OU}.  Its infinitesimal generator is
$L=\frac{s^2}{2}\frac{\d}{\d x^2}-\alpha x\frac{\d}{\d x}$, and the
eigenfunctions are known to be $H_k(x\sqrt{\alpha/s^2})$ where $H_k$ is
the Hermite polynomial $H_k(x):=(-1)^k \e^{x^2} \frac{\d^k}{\d x^k}
\e^{-x^2}$. The associated eigenvalue of $L$ is $-n\alpha$, so that the
spectral gap of $L$ is $\alpha$. Now the random walk we consider is the
flow $\e^{\deltat L}$ at time $\deltat$ of the process (with small
$\deltat$),
whose eigenvalues are $\e^{-n\alpha\deltat}$. So the spectral gap of the
discrete Laplacian $\e^{\deltat L}-\Id$ is $1-\e^{-\alpha\deltat}$. Since the
Ricci curvature is $1-\e^{-\alpha\deltat}$ too, the corollary is sharp again.

\section{Concentration results}

\subsection{Variance of Lipschitz functions}

We begin with the simplest kind of concentration, namely, an estimation
of the variance of Lipschitz functions. Contrary to Gaussian or
exponential concentration, the only assumption needed here is that the
average spread $\sigma$ is finite.

Since our Gaussian concentration result will yield basically the same
variance $\sigma^2/n\kappa$, we discuss sharpness of this estimate in
various examples in Section~\ref{sec:examplesrevisited}.

\begin{prop}
\label{prop:varlip}
Let $(X,d,m)$ be a random walk on a metric space, with Ricci curvature at least $\kappa>0$. Let $\nu$ be the unique invariant
distribution. Suppose that $\sigma<\infty$.

Then the variance of a $1$-Lipschitz function is at most
$\frac{\sigma^2}{n\kappa(2-\kappa)}\leq \frac{\sigma^2}{n\kappa}$.
\end{prop}

In particular, this implies that all Lipschitz functions are in
$L^2/\{\text{const}\}$; especially, $\int d(x,y)^2 \d\nu(x) \d\nu(y)$ is
finite. The fact that the Lipschitz norm controls the $L^2$ norm was used
above in the discussion of spectral properties of the random walk
operator.

\begin{dem}
Suppose for now that $f$ is bounded by $A\in \R$, so that $\Var f<\infty$.
We first prove that $\Var M^t\!f$ tends to $0$. Let $B_r$ be the ball of radius
$r$ in $X$ centered at some basepoint. Using that $M^t\!
f$ is $(1-\kappa)^t$-Lipschitz on $B_r$ and bounded by $A$ on
$X\!\,\setminus\!
B_r$, we get $\Var M^t\!f = \frac12 \iint (f(x)-f(y))^2\,\d\nu
(x)\d\nu(y)\leq
2(1-\kappa)^{2t}r^2+2A^2\nu(X\!\,\setminus\! B_r)$. Taking for example
$r=1/(1-\kappa)^{t/2}$ ensures that $\Var M^t\!f \to 0$.

As already mentioned, one has $\Var f = \Var M\!f + \int \Var_{m_x}
f\,\d\nu(x)$. Since $\Var M^t\!f \to 0$, by induction we get
\[
\Var f= \sum_{t=0}^\infty \int \Var_{m_x} M^t\!f \,\d\nu(x)
\]
Now by definition $\Var_{m_x} f \leq \sigma(x)^2/n_x$. Since $M^t\! f$ is
$(1-\kappa)^t$-Lipschitz, we have $\Var_{m_x} M^t\!f \leq
(1-\kappa)^{2t}\,
\sigma(x)^2/n_x$ so that the sum above is at most
$\frac{\sigma^2}{n\kappa(2-\kappa)}$. The case of unbounded $f$ is treated
by a simple limiting argument.
\end{dem}

\subsection{Gaussian concentration}

As mentioned above, positive Ricci curvature implies a
Gaussian-then-exponential concentration theorem. The estimated variance
is $\sigma^2/n\kappa$ as above, so that this is essentially a more
precise version of Proposition~\ref{prop:varlip}, with some loss in the
constants. We will see in the discussion below
(Section~\ref{sec:examplesrevisited}) that in the main
examples, the order of magnitude is correct.

The fact that concentration is not Gaussian far away from the mean is
genuine, as exemplified by the binomial distribution on the cube
(Section~\ref{sec:poisson}) or $M/M/\infty$ queues
(Section~\ref{sec:continuoustime}). A purely exponential behavior can be
achieved in very simple examples if $\sigma_\infty(x)$ is not bounded
(Example~\ref{ex:geom2}) or if the spread $\sigma(x)^2$ grows fast enough
(Section~\ref{sec:unstability}). In these examples, the transition from
Gaussian to non-Gaussian regime occurs roughly as predicted by the
theorem.

In the case of Riemannian manifolds, simply letting the step of
the random walk tend to $0$ makes the width of the Gaussian window tend
to infinity, so that we recover Gaussian concentration as in the
Lévy--Gromov or Bakry--Émery theorems.

The width of the Gaussian window is controlled by two factors: the
quantity $\sigma_\infty$, which represents the ``granularity'' of the
process and can result in Poisson-like behavior; and the rate of
variation of the spread $\sigma(x)^2$, which can result in exponential
behavior. The latter phenomenon yields to the assumption that
$\sigma(x)^2$ is bounded by a Lipschitz function.

\begin{thm}
\label{thm:gaussconc}
Let $(X,d,m)$ be an ergodic random walk on a metric space as above, with
invariant distribution $\nu$. Suppose that for any two distinct
points $x,y\in X$ one has $\kappa(x,y)\geq \kappa>0$.

Let
\[
D^2_x:=\frac{\sigma(x)^2}{n_x\kappa}
\]
and
\[
D^2:=\E_\nu D^2_x
\]

Suppose that the function $x\mapsto D^2_x$ is $C$-Lipschitz. Set
\[
t_{\text{max}}:=\frac{2D^2}{\max(2C,3\sigma_\infty)}
\]

Then for any $1$-Lipschitz function $f$, for any $t\leq t_{\text{max}}$
we have
\[
\nu\left(\left\{x,f(x)\geq t+\E_\nu f\right\}\right)\leq \exp
\,-\,\frac{t^2}{6D^2}
\]
and for $t\geq t_{\text{max}}$
\[
\nu\left(\left\{x,f(x)\geq t+\E_\nu f\right\}\right)\leq \exp
\,\left(-\,\frac{t_\text{max}^2}{6D^2}-\frac{t-t_\text{max}}{\max(2C,3\sigma_\infty)}\right)
\]
\end{thm}

\begin{rem}
\label{rem:finitetime}
It is clear from the proof below that $\sigma(x)^2/n_x\kappa$ itself need not be
Lipschitz, only bounded by some Lipschitz function. In particular, if
$\sigma(x)^2$ is bounded one can always take
$D^2=\sup_x\frac{\sigma(x)^2}{n_x\kappa}$ and $C=0$.
\end{rem}

It might seem that, in order to estimate $\E_\nu D_x^2$, one needs to
know in advance concentration properties for the invariant distribution
$\nu$; however, Proposition~\ref{prop:meandiam} or
Corollary~\ref{cor:avlip} often provides sharp estimates for $\E_\nu
D_x^2$, as we shall see in the examples.

In Secion~\ref{sec:unstability}, we give a simple example where the
Lipschitz constant of $\sigma(x)^2$ is large, resulting in exponential
rather than Gaussian behavior. In Section~\ref{sec:heavy} we give an
example of a process with quadratic growth of $\sigma(x)^2$, and which
exhibits non-exponential tails. Thus the Lipschitz assumption cannot
simply be removed.

The assumption that $\sigma_\infty$ is bounded can be replaced with a
Gaussian-type control for the local measures $m_x$, which however
generally results in much poorer estimates of the variance in discrete
situations (see Remark~\ref{rem:DGW}).

\begin{dem}
This proof is a variation on standard martingale methods for
concentration (see e.g.\ Lemma~4.1 in~\cite{Led01}).

Let $f$ be a $1$-Lipschitz function and $\lambda\geq 0$. For any
smooth
function $g$ and any real-valued random variable $Y$, a Taylor expansion
gives $\E g(Y)\leq g(\E
Y)+\frac12 (\sup g'') \Var Y$, so that
\[
(M\e^{\lambda f})(x)\leq \e^{\lambda
M\!f(x)}+\frac{\lambda^2 \e^{\lambda (M\!f(x)+2\sigma_\infty)}}{2}\Var_{m_x}
f
\]

Take $\lambda<1/3\sigma_\infty$ so that $\e^{2\lambda\sigma_\infty}\leq
2$. By definition, $\Var_{m_x} f\leq \norm{f}^2_{\text{Lip}}\sigma(x)^2/n_x$, hence
\[
(M\e^{\lambda f})(x)\leq \e^{\lambda
M\!f(x)}\left(1+\lambda^2\frac{\sigma(x)^2}{n_x}\right)\leq
\e^{\lambda \left(M\!f(x)+\lambda\frac{\sigma(x)^2}{n_x}\right)}
\]

But since $\sigma(x)^2/n_x\kappa$ is $C$-Lipschitz by assumption, and
since besides $M\!f(x)$ is $(1-\kappa)$-Lipschitz by
Proposition~\ref{prop:lipcontracting}, the sum $M\!f(x)+\lambda
\frac{\sigma(x)^2}{n_x}$ is $(1-\kappa+\lambda C \kappa)$-Lipschitz. 

From now on we take $\lambda\leq 1/2C$. We can repeat the argument,
setting $f_1(x):=M\!f(x)+\lambda \frac{\sigma(x)^2}{n_x}$ and using that
$f_1$ is $(1-\kappa/2)$-Lipschitz. This yields
\[
(M^2 \e^{\lambda f})(x)\leq (M\e^{\lambda f_1})(x)\leq
\e^{\lambda M\!f_1(x)+\lambda^2\frac{\sigma(x)^2}{n_x}(1-\kappa/2)^2}
\]

Next, $M\!f_1$ is $(1-\kappa)(1-\kappa/2)$-Lipschitz, whereas $\lambda
\frac{\sigma(x)^2}{n_x}(1-\kappa/2)^2$ is
$\frac{\kappa}{2}(1-\kappa/2)^2$-Lipschitz. So $f_2(x):=M\!f_1(x)+\lambda
\frac{\sigma(x)^2}{n_x}(1-\kappa/2)^2$ is (at least)
$(1-\kappa/2)^2$-Lipschitz, hence
\[
(M^3 \e^{\lambda f})(x)\leq (M\e^{\lambda f_2})(x)\leq
\e^{\lambda M\!f_2(x)+\lambda^2\frac{\sigma(x)^2}{n_x}(1-\kappa/2)^4}
\]

By induction, we get that
$f_{k+1}(x):=M\!f_k(x)+\lambda\frac{\sigma(x)^2}{n_x}(1-\kappa/2)^{2k}$
is $(1-\kappa/2)^{k+1}$-Lipschitz and that $(M^k\e^{\lambda f})(x)\leq
\e^{\lambda f_k(x)}$.

Now setting $g(x):=\frac{\sigma(x)^2}{n_x}$ and expanding $f_k$ yields
\[
f_k(x) = (M^k\!f)(x)+\lambda \sum_{i=1}^{k}
(M^{k-i}\!g)(x)\,(1-\kappa/2)^{2(i-1)}
\]
so that the limit of $f_k(x)$ when $k\to\infty$ is
\[
\E_\nu f+\lambda \sum_{i=1}^{\infty} \E_\nu g\,(1-\kappa/2)^{2(i-1)}
\leq \E_\nu f+\lambda \E_\nu g\,\frac{4}{3\kappa} 
\]

Meanwhile, $(M^k \e^{\lambda f})(x)$ tends to $\E_\nu \e^{\lambda f}$, so
that
\[
\E_\nu \e^{\lambda f}\leq \e^{\lambda \E_\nu
f+\frac{4\lambda^2}{3\kappa}\E_\nu \frac{\sigma(x)^2}{n_x}}
\]

We can conclude by a standard
Chebyshev inequality argument.
\end{dem}

\begin{rem}
The proof provides a similar concentration result for the finite-time measures
$\mu^{\ast k}_x$ as well, with variance
\[
D^2_{x,k}=\sum_{i=1}^k
(1-\kappa/2)^{2(i-1)}\left(M^{k-i}\frac{\sigma(y)^2}{n_y}\right)\!(x)
\]
and the same expression for $t_\text{max}$.
\end{rem}

\begin{rem}
\label{rem:DGW}
The condition that $\sigma_\infty$ is uniformly bounded can be replaced
with a Gaussian-type assumption, namely that for each measure $m_x$ there
exists a number $s_x$ such that $\E_{m_x}\e^{\lambda f}\leq \e^{\lambda^2
s_x^2/2} \e^{\lambda \E_{m_x}f}$ for any $1$-Lipschitz function $f$.  Then
a similar theorem holds, with $\sigma(x)^2$ replaced with $s_x^2$.
(When $s_x^2$ is constant this is Proposition~2.10 in~\cite{DGW04}.)
However, this is generally not well-suited to discrete settings, because
when transition probabilities are small, the best $s_x^2$ for which such
an
inequality is satisfied is usually much larger than the actual variance
$\sigma(x)^2$: for example, if two points $x$ and $y$ are at distance $1$
and $m_x(y)=\eps$, $s_x$ must satisfy $\e^{-1/2s_x^2}\leq \eps$ hence
$s_x^2\geq 1/2\ln(1/\eps)\gg \eps$. Thus making this assumption will
provide extremely poor estimates of the variance $D^2$ when some
transition probabilities are small (e.g.\ for binomial distributions on
the discrete cube); however, when this does not occur (e.g.\ for the
uniform distribution on the discrete cube), this assumption allows to get
rid of $\sigma_\infty$, and even get genuine Gaussian concentration for
all $t\in \R$ in the case $C=0$.
\end{rem}

\subsection{Examples revisited}
\label{sec:examplesrevisited}

Let us test the sharpness of these estimates in some examples, beginning
with the simplest ones. In each case, we gather the relevant quantities
in a table.  Recall that $\approx$ denotes an equality up to a
multiplicative universal constant (typically $\leq 4$), while symbol
$\sim$ denotes usual asymptotic equivalence (with the correct constant).

\newenvironment{extable}{\begin{center}\begin{tabular}{|ll|}\hline}{\\\hline\end{tabular}\end{center}}

\subsubsection{Riemannian manifolds}
\label{sec:riem}

First, let $X$ be a $N$-dimensional Riemannian manifold with positive
Ricci curvature.  Equip this manifold with the random walk at scale
$\eps>0$, as in Example~\ref{ex:riem}.

Let $\inf \Ric$ denote the largest $K>0$ such that $\Ric(v,v)\geq K$ for
any unit tangent vector $v$. The the relevant quantities for this random
walk are as follows (see Section~\ref{sec:riemproof} for the proofs).

\begin{extable}
Ricci curvature & $\kappa\sim \frac{\eps^2}{2(N+2)}\,\inf\Ric$
\\
Spread & $\sigma(x)^2\sim \eps^2 \frac{N}{N+2}\quad \forall x$
\\
Dimension & $n\approx N$
\\
Variance (Lévy--Gromov thm.) & $\approx 1/\inf\Ric$
\\
Gaussian variance (Thm.~\ref{thm:gaussconc}) & $D^2\approx 1/\inf \Ric$
\\
Gaussian range & $t_\text{max} \approx 1/(\eps\inf\Ric)\,\to\infty$
\end{extable}

So, up to some (small) constants, we recover Gaussian
concentration as in the Lévy-Gromov theorem.

The same applies to diffusions with a drift on a Riemannian manifold. To
be consistent with the notation of Example~\ref{ex:BE}, in the table
above $\eps$ has to be replaced with $\sqrt{(N+2) \deltat}$, and $\inf
\Ric$ with $\inf \left(\Ric(v,v)-2\nabla^{\text{sym}}F(v,v)\right)$ for
$v$ a unit tangent vector. (In the non-compact case, care has to be taken
since the Brownian motion on the manifold may not exist, and even if it
does its approximation at time $\deltat$ may not converge uniformly on
the manifold. In explicit examples such as the Ornstein--Uhlenbeck
process, however, this is not a problem.)

\subsubsection{Discrete cube}

Back to the discrete cube $\{0,1\}^N$ of Example~\ref{ex:cube}, equipped
with its graph distance (Hamming metric) and lazy random walk.

\begin{extable}
Ricci curvature & $\kappa=1/N$
\\
Spread & $\sigma(x)^2\approx 1\quad \forall x$
\\
Dimension & $n\approx 1$
\\
Gaussian variance (Thm.~\ref{thm:gaussconc}) & $D^2\approx N$
\\
Actual variance & $N/4$
\end{extable}

The following simple remark allows to actually compute the small
numerical constants implied in the notation $\approx$, and to check that
Proposition~\ref{prop:varlip} gives a sharp value when $N\to\infty$.

\begin{prop}
Let $m$ be the lazy simple random walk on a locally finite graph. Then,
for any vertex $x$ one has $\sigma(x)^2/n_x\leq 1/2$.
\end{prop}

Applying this to the estimate of Proposition~\ref{prop:varlip} for the
discrete cube, one gets $\sigma^2/n\kappa(2-\kappa)\leq
1/2\kappa(2-\kappa)$ which, for $\kappa=1/N$, yields $N/2(2-1/N)\sim
N/4$. (One can actually get exactly $N/4$ by using a continuous-time
random walk instead.)

\begin{dem}
By definition $\sigma(x)^2/n_x$ is the maximal variance, under $m_x$, of a
$1$-Lipschitz function. So let $f$ be a $1$-Lipschitz function on the
graph. Since variance is unvariant by adding a constant, we can assume
that $f(x)=0$. Then $\abs{f(y)}\leq 1$ for any neighbor $y$ of $x$. Since
$m$ is the lazy simple random walk, we have $m_x(x)\geq 1/2$ (with
equality if there are no loops) and the mass, under $m_x$, of all
neighbors of $x$ is at most $1/2$. Hence $\Var_{m_x}
f=\E_{m_x}f^2-(\E_{m_x}f)^2\leq \E_{m_x}f^2\leq 1/2$.

This value is actually achieved when $x$ has an even number of neighbors
and when no two distinct neighbors of $x$ are neighbors; in this case one
can take $f(x)=0$, $f=1$ on half the neighbors of $x$ and $f=-1$ on the
remaining neighbors of $x$.
\end{dem}

%
%
%

\subsubsection{Binomial distributions}
\label{sec:poisson}

The occurrence of a finite range $t_\text{max}$ for the Gaussian behavior
of tails is genuine, as the following example shows.

Let $X=\{0,1\}^N$ equipped with its Hamming metric (each edge is of
length $1$). Consider the following Markov chain on $X$: for some
$0<p<1$, at each step, choose a bit at random among the $N$ bits; if it
is equal to $0$, flip it to $1$ with probability $p$; if it is equal to
$1$, flip it to $0$ with probability $1-p$. The binomial distribution
$\nu\left((x_i)\right)=\prod p^{x_i}(1-p)^{1-x_i}$ is reversible for this Markov
chain. The Ricci curvature of this Markov chain is $1/N$.

Let $k$ be the number of bits of $x\in X$ which are equal to $1$. Then
$k$ follows a Markov chain on $\{0,1,\ldots,N\}$, whose transition probabilities are:
\begin{align*}
p_{k,k+1}&=p(1-k/N)
\\
p_{k,k-1}&=(1-p)k/N
\\
p_{k,k}&=pk/N+(1-p)(1-k/N)
\end{align*}

The binomial distribution with parameters $N$ and $p$, namely
$\binom{N}{k} p^k(1-p)^{N-k}$, is reversible for
this Markov chain. Moreover, the Ricci curvature of this Markov chain is
$1/N$.

Now, fix some $\lambda>0$ and consider the case $p=\lambda/N$. Let
$N\to\infty$. It is well-known that the invariant distribution tends to
the Poisson distribution $\e^{-\lambda}\lambda^k/k!$ on $\N$. 

Let us see how Theorem~\ref{thm:gaussconc} performs on this example. The
table below applies either to the full space $\{0,1\}^N$, with $k$ the
function ``number of $1$'s'', or to its projection on $\{0,1,\ldots,N\}$.
Note the use of Proposition~\ref{prop:meandiam} to estimate $\sigma^2$,
without having to resort to explicit knowledge of the invariant distribution.
(All constants implied in the $O(1/N)$ notation are small and completely
explicit.)  

\begin{extable}
Ricci curvature & $\kappa=1/N$
\\
Spread & $\sigma(k)^2=(\lambda+k)/N+O(1/N^2)$
\\
Estimated  $\E k$ (Prop.~\ref{prop:meandiam})& $\E k\leq
J(0)/\kappa=\lambda$
\\
Actual $\E k$ & $\E k=\lambda$
\\
Average spread & $\sigma^2=\E \sigma(k)^2= 2\lambda/N+O(1/N^2)$
\\
Dimension & $n\geq 1$
\\
Estimated variance (Prop.~\ref{prop:varlip}) &
$\sigma^2/n\kappa(2-\kappa)= \lambda+O(1/N)$
\\Actual variance & $\lambda$
\\Gaussian variance (Thm.~\ref{thm:gaussconc}) & $D^2= 2\lambda+O(1/N)$
\\Lipschitz constant of $D_x^2$ & $C= 1+O(1/N)$
\\Gaussian range & $t_\text{max}=4\lambda/3$
\end{extable}

The Poisson distribution has a roughly Gaussian behavior (with variance
$\lambda$) in a range of size approximately $\lambda$ around the mean;
further away, it decreases like $\e^{-k\ln k}$ which is not Gaussian. This
is in good accordance with the theorem, and shows that the Gaussian range
cannot be extended.

\subsubsection{A continuous-time example: $M/M/\infty$ queues}
\label{sec:continuoustime}

Here we show how to apply the theorem above to a continuous-time example,
the $M/M/\infty$ queue. These queues were brought to my attention by
D.~Chafaï.

The $M/M/\infty$ queue consists in an infinite number of ``servers''. Each
server can be free ($0$) or busy ($1$). The state space consists in all
sequences in $\{0,1\}^\N$ with a finite number of $1$'s. The dynamics is
at follows: Fix two numbers $\lambda>0$ and $\mu>0$. At a rate $\lambda$
per unit of time, a client arrives and the first free server becomes
busy. At a rate $\mu$ per unit of time, each busy server finishes its job
(independently of the others) and becomes free. The number $k\in\N$ of
busy servers is a continuous-time Markov chain, whose transition
probabilities at small times $t$ are given by
\begin{align*}
p^t_{k,k+1}&=\lambda t+O(t^2)
\\
p^t_{k,k-1}&=k\mu t+O(t^2)
\\
p^t_{k,k}&=1-(\lambda +k\mu) t + O(t^2)
\end{align*}

If we replace $\lambda$ with $\lambda/N$ and $\mu$ with $1/N$, this
Markov chain appears as the limit of the binomial example above. This is
especially clear in the table below.

This system is often presented as a discrete analogue of an
Ornstein--Uhlenbeck process, since asymptotically the drift is linear
towards the origin. However, it is not symmetric around the mean, and
moreover the invariant (actually reversible) distribution $\nu$ is a Poisson
distribution (with parameter $\lambda/\mu$), rather than a Gaussian.

In this continuous-time setting, the definition are adapted as follows:
$\kappa(x,y):=-\,\frac{\d}{\d t} \Td(m^t_x,m^t_y)/d(x,y)$ (as mentioned
in the introduction) and $\sigma(x)^2:=\frac12 \frac{\d}{\d t}\iint
d(y,z)\,\d m^t_x(y) \d m^t_x(z)$, where $m^t_x$ is the law at time $t$ of
the process starting at $x$. It is immediate to check that
the Ricci curvature of this process is $\mu$.
Proposition~\ref{prop:varlip} (with
$\sigma^2/2n\kappa$ instead of $\sigma^2/n\kappa(2-\kappa)$ because both
$\sigma^2$ and $\kappa$ tend to $0$ for the discrete-time approximation)
and Theorem~\ref{thm:gaussconc} still hold. 

The relevant quantities are as follows.

\begin{extable}
Ricci curvature & $\kappa=\mu$
\\
Spread & $\sigma(k)^2=k\mu+\lambda$
\\
Estimated  $\E k$ (Prop.~\ref{prop:meandiam})& $\E k\leq
J(0)/\kappa=\lambda/\mu$
\\
Actual $\E k$ & $\E k=\lambda/\mu$
\\
Average spread & $\sigma^2=\E \sigma(k)^2= 2\lambda$
\\
Dimension & $n\geq 1$
\\
Estimated variance (Prop.~\ref{prop:varlip}) &
$\sigma^2/2n\kappa= \lambda/\mu$
\\Actual variance & $\lambda/\mu$
\\
Gaussian variance (Thm.~\ref{thm:gaussconc}) & $D^2= 2\lambda/\mu$
\\
Lipschitz constant of $D_x^2$ & $C=1$
\\
Gaussian range & $t_\text{max}=4\lambda/3\mu$
\end{extable}

So once more Theorem~\ref{thm:gaussconc} is in excellent accordance with
the behavior of the random walk, whose invariant distribution is Poisson
with mean and variance $\lambda/\mu$.

An advantage of this approach is that is can be generalized to situations
where the rates of the servers are not constant, but, say, bounded
between, say, $\mu_0/10$ and $10\mu_0$. Indeed, the $M/M/\infty$ queue
above can be seen as a Markov chain in the full configuration space of
the servers, namely the space of all sequences over the alphabet
$\{\text{free},\text{busy}\}$ containing a finite number of ``busy''. It
is easy to check that the Ricci curvature is still equal to $\mu$ in this
configuration space. Now let us consider the case of variable rates: in
this situation, the number of busy servers is generally not Markovian, so
one has to work in the configuration space. If the rate of the $i$-th
server is $\mu_i$, the Ricci curvature is $\inf \mu_i$ in the
configuration space, whereas the spread is controlled by $\sup \mu_i$. So
if the rates vary in a bounded range, Ricci curvature still provides a
Gaussian-type control, though an explicit description of the invariant
distribution is not available.

\subsubsection{An example of exponential concentration}
\label{sec:unstability}

We give here a very simple example of a Markov chain which has positive
curvature but for which concentration is not Gaussian but exponential,
due to large variations of the spread, resulting in a large value of $C$.
An even simpler example, with exponential concentration due to unbounded
$\sigma_\infty(x)$, was given in the introduction
(Example~\ref{ex:geom2}).

This is a continuous-time random walk on $\N$ defined as follows. Take
$\alpha<\beta\in\R$.  For $k\in \N$, the transition rate from $k$ to
$k+1$ is $(k+1)\alpha$, whereas the transition rate from $k+1$ to $k$
is $(k+1)\beta$. It is immediate to check that the geometric distribution
with decay $\alpha/\beta$ is reversible for this Markov chain.

The Ricci curvature of this Markov chain is easily seen to be
$\beta-\alpha$. We have $\sigma(k)^2=(k+1)\alpha+k\beta$, so that
$\sigma(k)^2$ is $(\alpha+\beta)$-Lipschitz and
$C=(\alpha+\beta)/(\beta-\alpha)$.

The expectation of $k$ under the invariant distribution can be bounded by
$J(0)/\kappa=\alpha/(\beta-\alpha)$ by Proposition~\ref{prop:meandiam},
which is actually the exact value. So the expression above for
$\sigma(k)^2$ yields $\sigma^2=2\alpha\beta/(\beta-\alpha)$.
Consequently, the estimated variance $\sigma^2/2n\kappa$ (obtained by the
continuous-time version of Proposition~\ref{prop:varlip}) is at most
$\alpha\beta/(\beta-\alpha)^2$, which is the actual value.

Now consider the case when $\beta-\alpha$ is small.  If we try to apply
Theorem~\ref{thm:gaussconc} without taking into account the variations of
the spread (witnessed by the constant $C$), we get blatantly false
results since the invariant distribution is not Gaussian at all. In the
regime where $\beta-\alpha\to 0$, the width of the Gaussian window in
Theorem~\ref{thm:gaussconc} is $D^2/C\approx\alpha/(\beta-\alpha)$. This
is fine, as this is the decay distance of the invariant distribution, and
in this interval both the Gaussian and geometric estimates are close to
$1$ anyway.
But if the $C$ factor was not included, we would get
$D^2/\sigma_\infty=\alpha\beta/(\beta-\alpha)^2$, which is much larger;
the invariant distribution is clearly not Gaussian on this interval.

Moreover, Theorem~\ref{thm:gaussconc} predicts, in the exponential
regime, a $\exp(-t/2C)$ behavior for concentration. Here the asymptotic
behavior of the invariant distribution is
$(\alpha/\beta)^t\sim(1-2/C)^t\sim \e^{-2t/C}$ when $\beta-\alpha$ is
small. So we see that (up to a constant $4$) the exponential decay rate
predicted by Theorem~\ref{thm:gaussconc} is genuine.

\subsubsection{Heavy tails}
\label{sec:heavy}

It is clear that a variance control alone does not imply any
concentration beyond the Bienaymé-Chebyshev inequality. We now show that
this is till the case even with the positive curvature assumption. Namely,
in Theorem~\ref{thm:gaussconc}, neither the assumption that $\sigma(x)^2$
is Lipschitz, nor the assumption that $\sigma_\infty$ is bounded, can be
removed (but see Remark~\ref{rem:DGW}).

\paragraph{Heavy tails with non-Lipschitz $\sigma(x)^2$.}
Our next example shows that if the spread $\sigma(x)^2$ is not Lipschitz,
then non-exponential tails may occur in spite of positive curvature.

Consider the continuous-time random walk on $\N$ defined as follows: the
transition rate from $k$ to $k+1$ is $a(k+1)^2$, whereas the transition
rate from $k$ to $k-1$ is $a(k+1)^2+bk$ for $k\geq 1$. Here $a,b>0$ are
fixed.

We have $\kappa=b$ and $\sigma(k)^2=2a(k+1)^2+bk$, which is obviously not
Lipschitz. 

This Markov chain has a reversible measure $\nu$, which satisfies
$\nu(k)/\nu(k-1)=ak^2/(a(k+1)^2+bk)=1-\frac1k(2+\frac{b}{a})+O(1/k^2)$.
Consequently, asymptotically $\nu(k)$ behaves like \[
\prod_{i=1}^k\left(1-\tfrac1i(2+\tfrac{b}{a})\right)\approx
\e^{-(2+b/a)\sum_{i=1}^k\frac1i}\approx k^{-(2+b/a)}
\]
thus exhibiting heavy, non-exponential tails.

This shows that the Lipschitz assumption for $\sigma(x)^2$ cannot be
removed, even if in this case $\sigma_\infty$ is bounded by $1$. It would seem reasonable to look for a systematic correspondance
between the asymptotic behavior of $\sigma(x)^2$ and the behavior of
tails.

\paragraph{Heavy tails with unbounded $\sigma_\infty$.} Consider the
following random walk on $\N^\ast$: a number $k$ goes to $1$ with
probability $1-1/4k^2$ and to $2k$ with probability $1/4k^2$. One can
check that $\kappa\geq 1/2$. These probabilities are chosen so that
$\sigma(k)^2=(2k-1)^2\times 1/4k^2\times (1-1/4k^2)\leq 1$, so that the
variance of the invariant distribution is small. However, let us evaluate
the probability that, starting at $1$, the first $i$ steps consist in
doing a multiplication by $2$, so that we end at $2^i$; this probability
is $\prod_{j=0}^{i-1}\frac1{4.(2^j)^2}=4^{-1-i(i-1)/2}$. Setting
$i=\log_2 k$, we see that the invariant distribution $\nu$ satisfies
\[
\nu(k)\geq \frac{\nu(1)}{4}\,2^{-\log_2 \!k \,(\log_2 \!k-1)}
\]
for $k$ a power of $2$.
This is clearly not Gaussian or exponential, though $\sigma(k)^2$ is
bounded.

\section{Local control and logarithmic Sobolev inequality}

The estimates above (e.g.\ for the spectral gap) were global: we used
that the averaging operator $M$ transforms a $1$-Lipschitz function into a
$(1-\kappa)$-Lipschitz function. Now we turn to some form of control of
the gradient of $M\!f$ at some point, in terms of the gradient of $f$ at
neighboring points. This is closer to classical Bakry--Émery theory, and
allows to get a kind of logarithmic Sobolev inequality.

\begin{defi}
\label{def:gradient}
Choose $\lambda>0$ and, for any function $f:X\to \R$, define the
\emph{$\lambda$-range gradient} of $f$ by
\[
(Df)(x):=\sup_{y,y'\in X}\frac{\abs{f(y)-f(y')}}{d(y,y')}\,\e^{-\lambda
d(x,y)-\lambda d(x,y')}
\]
\end{defi}

This is a kind of ``mesoscopic'' Lipschitz constant of $f$ around $x$.
Note that if $f$ is a smooth function on a compact Riemannian manifold,
when $\lambda\to \infty$ this quantity tends to $\abs{\nabla \!f(x)}$.

It is important to note that $Df$ is $2\lambda$-log-Lipschitz.

We will also need a control on negative curvature: In a Riemannian
manifold, the Ricci curvature might be $\geq \eps$ because there is a
direction of curvature $1$ and a direction of curvature $-1+\eps$. The
next definition captures these variations.

\begin{defi}[ (Unstability)]
\label{def:unstability}
Let
\[
\kappa_+(x,y):=\frac{1}{d(x,y)}\int_z (d(x,y)-d(x+z,y+z))_+
\]
and
\[
\kappa_-(x,y):=\frac{1}{d(x,y)}\int_z (d(x,y)-d(x+z,y+z))_-
\]
where $a_+$ and $a_-$ are the positive and negative part of $a\in \R$, so
that $\kappa(x,y)=\kappa_+(x,y)-\kappa_-(x,y)$. (The integration over $z$
is under a coupling realizing the value of $\kappa(x,y)$.)

The unstability $U(x,y)$ is defined as
\[
U(x,y):=\frac{\kappa_-(x,y)}{\kappa(x,y)}\qquad\text{and}\qquad
U:=\sup_{x,y \in X, \,x\neq y} U(x,y)
\]
\end{defi}


\begin{rem}
If $X$ is $\eps$-geodesic, then an upper bound for $U(x,y)$ with
$d(x,y)\leq \eps$ implies the same upper bound for $U$.
\end{rem}

In most discrete examples given in the introduction
(Examples~\ref{ex:cube}, \ref{ex:discrOU}, \ref{ex:mult}, \ref{ex:geom},
\ref{ex:geom2}), unstability is actually $0$, meaning that the coupling
between $m_x$ and $m_y$ never increases distances (this could be a
possible definition of non-negative sectional curvature for Markov
chains).  In Riemannian manifolds, unstability is controlled by the
largest negative sectional curvature, but this does not influence the
final results since one can take arbitrarily small steps for the random
walk.  Interestingly, in Example~\ref{ex:glauber} (Glauber dynamics),
unstability depends on temperature.

Due to the use of the gradient $D$, the theorem below is interesting only
if a reasonable estimate for $Df$ can be obtained depending on ``local''
data. This is not the case when $f$ is not $\lambda$-log-Lipschitz. This
is consistent with the fact mentioned above, that Gaussian concentration
of measure only occurs in a finite range, with exponential concentration
afterwards, which implies that no true logarithmic Sobolev inequality can
hold in general.

\begin{thm}
\label{thm:logsob}
Suppose that Ricci curvature is at least $\kappa>0$.
Let $\lambda\leq\frac{1}{24\sigma_\infty(1+U)}$ and consider the $\lambda$-range
gradient $Df$. Then for any function $f:x\to \R$ such that
$Df<\infty$, one has
\[
\Var_\nu f\leq \left(\sup_x \frac{4\sigma(x)^2}{\kappa n_x}\right) \int
(Df)^2 \,\d\nu
\]
and for positive $f$,
\[
\Ent_\nu f\leq \left(\sup_x \frac{4\sigma(x)^2}{\kappa n_x}\right) \int
\frac{(Df)^2}{f} \,\d\nu
\]
where $\nu$ is the invariant distribution.

If moreover the random walk is reversible with respect to $\nu$, then
\[
\Var_\nu f\leq \int V(x) \,Df(x)^2\,\d\nu(x)
\]
and
\[
\Ent_\nu f \leq \int V(x) \,\frac{Df(x)^2}{f(x)}\,\d\nu(x)
\]
where
\[
V(x)=2\sum_{t=0}^\infty (1-\kappa/2)^{2t}
M^{t+1}\!\!\left(\frac{\sigma(x)^2}{n_x}\right)
\]
\end{thm}

The form involving $V(x)$ is motivated by the fact that, for reversible
diffusions in $\R^N$ with non-constant diffusion coefficients, these
coefficients naturally appear in the formulation of functional
inequalities (see e.g.~\cite{AMTU01}). The quantity $V(x) \,Df(x)^2$ is
to be thought of as a crude version of the Dirichlet form associated with
the random walk. It would be more satisfying to obtain inequalities
involving the latter (compare Corollary~\ref{cor:poincare}), but I could
not get a version of the commutation property $DM\leq (1-\kappa/2)MD$
involving the Dirichlet form.

\begin{rem}
\label{rem:V(x)}
If $\frac{\sigma(x)^2}{n_x\kappa}$ is $C$-Lipschitz (as in
Theorem~\ref{thm:gaussconc}),
then 
$V(x)\leq \frac{4}{\kappa}\int
\frac{\sigma(x)^2}{n_x}\,\d\nu(x)+2C\frac{J(x)}{\kappa}$.
\end{rem}

\paragraph{Examples.} Let us compare this theorem to classical results.

In the case of a Riemannian manifold, for any smooth function $f$ we can
choose a random walk with small enough steps, so that $\lambda$ can be
arbitrarily large and $Df$ arbitrarily close to $\abs{\nabla\! f}$. Since
moreover $\sigma(x)^2$ does not depend on $x$ for the Brownian motion, this
theorem allows to recover the logarithmic Sobolev inequality in the
Bakry--Émery framework, with the correct constant up to a factor $4$.

Now consider the two-point space $\{0,1\}$, equipped with the measure
$\nu(0)=1-p$ and $\nu(1)=p$. This is a classical space on which modified
logarithmic Sobolev inequalities were introduced~\cite{BL98}. We endow
this space with the Markov chain sending each point to the invariant
distribution. Here we have $\sigma(x)^2=p(1-p)$, $n_x=1$ and $\kappa=1$,
so that we get the inequality $\Ent_\nu f\leq 4p(1-p)\int
\frac{(Df)^2}{f}\,\d\nu$, identical to the known inequality~\cite{BL98} except
for the factor $4$.

Tensorizing this result provides a modified logarithmic inequality for
Bernoulli and Poisson measures~\cite{BL98}. If, instead, we directly
apply the theorem above to the Bernoulli measure on $\{0,1\}^N$ or the
Poisson measure on $\N$ (see Sections~\ref{sec:poisson}
and~\ref{sec:continuoustime}), we get slightly worse results. Indeed,
consider the $M/M/\infty$ queue on $\N$, which is the limit when
$N\to\infty$ of the projection on $\N$ of the Markov chains on
$\{0,1\}^N$ associated with Bernoulli measures.
Keeping the notation of Section~\ref{sec:continuoustime}, we get, in the
continuous-time version, $\sigma(x)^2=x\mu+\lambda$, which is not
constant. So we have to use $V(x)$; Remark~\ref{rem:V(x)} and the
formulas in Section~\ref{sec:continuoustime} yields
$V(x)\leq 8\lambda/\mu+2(\lambda+x\mu)/\mu$ so that we get the inequality
\[
\Ent_\nu f\leq \frac{\lambda}{\mu}\int
\frac{Df(x)^2}{f(x)}\,(10+2x\mu/\lambda)\,\d\nu(x)
\]
which is to be compared to the inequality
\[
\Ent_\nu f\leq \frac{\lambda}{\mu} \int \frac{D_{\!+}f(x)^2}{f(x)} \,\d\nu(x)
\]
obtained in~\cite{BL98}, with $D_{\!+}f(x)=f(x+1)-f(x)$. So
asymptotically our version is worse by a factor $x$. Note however that
the Poisson measure satisfies $x\mu/\lambda\,\,\d\nu(x)=\d\nu(x-1)$, so
one could say that our general, non-local notion of gradient fails to distinguish
between a point and an immediate neighbor, and does not take advantage of
the particular structure of a random walk on $\N$.

\paragraph{Proof.} We now turn to the proof of Theorem~\ref{thm:logsob},
which is essentially a copy of the Bakry--Émery argument. The key
property is Proposition~\ref{prop:localcontrol}, a commutation property
between the gradient and random walk operators stating that $DM\leq
(1-\kappa/2)MD$.

\begin{lem}
\label{lem:localdecr}
Let $A$ be a function on $\Supp m_x$, such that
$A(z)\leq \e^\rho A(z')$ for any $z,z'\in \Supp m_x$, with
$\rho\leq 
\frac{1}{2(1+U)}$. Then for any $x,y\in X$ we have
\[
\int_z A(z)\frac{d(x+z,y+z)}{d(x,y)}\leq (1-\kappa(x,y)/2)\int_z A(z)
\]
and in particular
\[
\int_z A(z)(d(x+z,y+z)-d(x,y))\leq 0
\]
\end{lem}

\begin{dem}
Set $F=\max_z A(z)$. Then
\begin{align*}
\int_z A(z)\frac{d(x+z,y+z)}{d(x,y)}
&=\int_z A(z)+F\int_z
\frac{A(z)}{F}\left(\frac{d(x+z,y+z)}{d(x,y)}-1\right)
\end{align*}
and recall that, by definition, $\kappa_-(x,y)=\int_{z, d(x+z,y+z)>d(x,y)}
\left(d(x+z,y+z)/d(x,y)-1\right)$ and $\kappa_+(x,y)=\int_{z,
d(x+z,y+z)\leq d(x,y)}\left(1-d(x+z,y+z)/d(x,y)\right)$. Using that
$A(z)\leq F$ on one hand and $A(z)\geq \e^{-\rho}F$ on the
other hand, we get
\[
\int_z A(z)\frac{d(x+z,y+z)}{d(x,y)}\leq
\int_z A(z)+F(\kappa_-(x,y) -\e^{-\rho}\kappa_+(x,y))
\]

Now, recall that by definition of $U$ we have $\kappa_-(x,y)\leq U
\kappa(x,y)$. It is not difficult to check that $\rho\leq
\frac{1}{2(1+U)}$
is enough to ensure that $\e^{-\rho}\kappa_+(x,y)-\kappa_-(x,y)\geq \kappa(x,y)/2$, hence 
\begin{align*}
\int_z A(z)\frac{d(x+z,y+z)}{d(x,y)}
&\leq
\int_z A(z)-F\kappa(x,y)/2
\\&\leq \int_z A(z)\,(1-\kappa(x,y)/2)
\end{align*}
as needed.
\end{dem}

\begin{prop}
\label{prop:localcontrol}
Suppose that the Ricci curvature is at least $\kappa>0$, and choose some
$\lambda\leq \frac{1}{24\sigma_\infty(1+U)}$. Then for any function $f:X\to
\R$ we have
\[
D(M\!f)(x)\leq (1-\kappa/2)M(Df)(x)
\]
\end{prop}

\begin{dem}
For any $y,y'\in X$ we have
\begin{align*}
&\frac{\abs{M\!f(y)-M\!f(y')}}{d(y,y')}\,\e^{-\lambda(d(x,y)+d(x,y'))}
\\&\leq
\int_z \abs{f(y+z)-f(y'+z)}\,\frac{\e^{-\lambda(d(x,y)+d(x,y'))}}{d(y,y')}
\\&\leq
\int_z
Df(x+z)\frac{d(y+z,y'+z)}{\e^{-\lambda(d(x+z,y+z)+d(x+z,y'+z))}}\,\frac{\e^{-\lambda(d(x,y)+d(x,y'))}}{d(y,y')}
\\&=
\int_z A(z)B(z)\,\frac{d(y+z,y'+z)}{d(y,y')}
\end{align*}
where
$A(z)=Df(x+z)$ and
$B(z)=\e^{\lambda(d(x+z,y+z)-d(x,y)+d(x+z,y'+z)-d(x,y'))}$.

For any $z$ we have $(1-\kappa(x,y))d(x,y)-4\sigma_\infty\leq
d(x+z,y+z)\leq (1-\kappa(x,y))d(x,y)+4\sigma_\infty$ and likewise for $y'$, so
that $B$ varies by a factor at most $\e^{8\lambda \sigma_\infty}$. Likewise,
since $Df$ is $2\lambda$-log-Lipschitz, $A$ varies by a factor at most
$\e^{4\lambda \sigma_\infty}$. So the quantity $A(z)B(z)$ varies by at most
$\e^{12\lambda \sigma_\infty}$.

So if $\lambda\leq \frac{1}{24\sigma_\infty(1+U)}$, we can apply
Lemma~\ref{lem:localdecr} and get
\[
\int_z A(z)B(z)\,\frac{d(y+z,y'+z)}{d(y,y')}\leq (1-\kappa/2) \int_z
A(z)B(z)
\]

Now we have $\int_z A(z)B(z)=\int_z A(z)+\int_z A(z)(B(z)-1)$. Unwinding
$B(z)$ and using that $\e^a-1\leq ae^a$ for any $a\in \R$, we get
\begin{align*}
&\int_z A(z)(B(z)-1)\leq
\\&
\lambda \int_z A(z)B(z)\left(d(x+z,y+z)-d(x,y)+d(x+z,y'+z)-d(x,y')\right)
\end{align*}
which is non-positive by Lemma~\ref{lem:localdecr}. Hence $\int_z
A(z)B(z)\leq \int_z A(z)$, which ends the proof.
\end{dem}

Let
$\nu$ be the invariant distribution.
Let $f$ be a positive function with $\int f\,\d\nu=1$. We know that
\begin{align*}
\Ent f&=\int_x M\!f(x)\left(\Ent_{m_x}\frac{f}{M\!f(x)}\right) \d\nu(x)+\Ent M\!f
\\&=
\sum_{t\geq 0}
\int_x M^{t+1}\!f(x)\left(\Ent_{m_x} \frac{M^t\!f}{M^{t+1}\!f(x)}\right) \d\nu(x)
\end{align*}
and similarly
\[
\Var f=\sum_{t\geq 0} \int_x \Var_{m_x}M^t\!f \,\d\nu(x)
\]

Now for any $y,z\in\Supp m_x$ we have $\abs{f(y)-f(z)}\leq
Df(y)d(y,z)\e^{\lambda d(y,z)}$. Since $Df$ is $2\lambda$-log-Lipschitz,
we have $Df(y)\leq \e^{4\lambda \sigma_\infty} M(Df)(x)$, so that
$\abs{f(y)-f(z)}\leq d(y,z)\,M(Df)(x)\,\e^{6\lambda \sigma_\infty}$, i.e.\ $f$
is $M(Df)(x)\,\e^{6\lambda \sigma_\infty}$-Lipschitz. Consequently
\[
\Var_{m_x} f\leq \frac{2 (M(Df)(x))^2 \,\sigma(x)^2}{n_x}
\]
and, using that $a\log a\leq a^2-a$, we get that
$\Ent_{m_x}\frac{f}{M\!f(x)}\leq \frac{1}{M\!f(x)^2}\Var_{m_x} f$
so
\[
\Ent_{m_x} \frac{f}{M\!f(x)}\leq \frac{2 (M(Df)(x))^2
\,\sigma(x)^2}{n_x\,M\!f(x)^2}
\]

Thus
\[
\Var f\leq 2\sum_{t\geq 0} \int_x \frac{\sigma(x)^2}{n_x}\,(M(DM^t\!f)(x))^2
\,\d\nu(x)
\]
and
\begin{align*}
\Ent f&\leq
2\sum_{t\geq 0}
\int_x \frac{\sigma(x)^2}{n_x}\,\frac{(M(DM^tf)(x))^2}{M^{t+1}f(x)}
\,\d\nu(x)
\end{align*}

By Proposition~\ref{prop:localcontrol}, we have $(DM^t\!f)(y)\leq (1-\kappa/2)^t
M^t(Df)(y)$, so that
\[
\Var f\leq 2\sum_{t\geq 0}
\int_x \frac{\sigma(x)^2}{n_x}\,
(M^{t+1}Df(x))^2\,(1-\kappa/2)^{2t}
\,\d\nu(x)
\]
and
\[
\Ent f\leq 2\sum_{t\geq 0}
\int_x \frac{\sigma(x)^2}{n_x}\,
\frac{(M^{t+1}Df(x))^2}{M^{t+1}\!f(x)} \,(1-\kappa/2)^{2t}
\,\d\nu(x)
\]

Now since the norm of $M$ acting on $L^2(\nu)$ is at most $1$, we have
\begin{align*}
\Var f&\leq 2\sup_x \frac{\sigma(x)^2}{n_x}\,\sum_{t\geq 0}(1-\kappa/2)^{2t}\int_x
(M^{t+1}Df(x))^2 \,\d\nu(x)
\\&\leq 
\frac{4}{\kappa} \sup_x\frac{\sigma(x)^2}{n_x}
\,\int_x (Df(x))^2 \,\d\nu(x)
\end{align*}

For the entropy of $f$, the Cauchy--Schwarz inequality yields
\[
(M^{t+1}Df(x))^2=\left(M^{t+1}\!\left(\frac{Df}{\sqrt{f}}.\sqrt{f}\right)\!(x)\right)^2
\leq M^{t+1}\!\left(\frac{(Df)^2}{f}\right)\!(x) \,M^{t+1}\!f(x)
\]
so that finally
\begin{align*}
\Ent f&\leq 2\sum_{t\geq 0}
\int_x \frac{\sigma(x)^2}{n_x}\,
M^{t+1}\!\left(\frac{(Df)^2}{f}\right)\!(x) \,(1-\kappa/2)^{2t}
\,\d\nu(x)
\\&\leq
\frac{4}{\kappa} \sup_x\frac{\sigma(x)^2}{n_x}
\,\int_x \frac{(Df(x))^2}{f(x)} \,\d\nu(x)
\end{align*}

%
%
%
%

\section{Exponential concentration in non-negative curvature}
\label{sec:exp}

We have seen that positive Ricci curvature implies a kind of Gaussian
concentration. We now show that non-negative Ricci curvature and the
existence of an ``attracting point'' imply exponential concentration.

The basic example to keep in mind is the following. Let $\N$ be the set
of non-negative integers equipped with its standard distance. Let $0<p<1$
and let the nearest-neighbor random walk on $\N$ that goes to the left
with probability $p$; explicitly $m_k=p\delta_{k-1}+(1-p)\delta_{k+1}$
for $k\geq 1$, and $m_0=p\delta_0+(1-p)\delta_1$.

Since for $k\geq 1$ the transition kernel is translation-invariant, it
is immediate to check that $\kappa(k,k+1)=0$. Besides, $\kappa(0,1)=p$.
There exists a invariant distribution if and only if $p>1/2$, and
it satisfies exponential concentration with characteristic decay distance
$1/\log(p/(1-p))$. For $p=1/2+\eps$ with small $\eps$ this behaves like
$1/4\eps$.

Geometrically, what entails exponential concentration in this example is
the fact that, for $p>1/2$, the point $0$ ``pulls'' its neighbor, and the
pulling is transmitted by non-negative Ricci curvature. We now formalize
this situation in the following theorem.

\begin{thm}
\label{thm:expconc}
Let $(X,d,(m_x))$ be a metric space with random walk. Suppose that
for some $o\in X$ and $r>0$ one has:
\begin{itemize}
\item $\kappa(x,y)\geq 0$ for all $x,y\in X$,
\item for all $x\in X$ with $r\leq d(o,x)<2r$, one has
$\Td(m_x,\delta_o)<d(x,o)$,
\item $X$ is $r$-geodesic,
\item There exists $s>0$ such that each measure $m_x$ satisfies the
Gaussian-type Laplace transform inequality
\[
\int \e^{\lambda f}\,\d m_x\leq \e^{\lambda^2s^2/2}\e^{\lambda \int f \d m_x}
\]
for any $\lambda>0$ and any $1$-Lipschitz function $f:\Supp m_x\to \R$.
\end{itemize}

Set $\rho=\inf\{d(x,o)-\Td(m_x,\delta_o), \,r\leq d(o,x)<2r\}$ and assume
$\rho>0$.

Then there exists a invariant distribution for the random walk.
Moreover, setting
$D=s^2/\rho$ and $m=r+2s^2/\rho+\rho(1+J(o)^2/4s^2)$, 
for any invariant distribution $\nu$
we have
\[
\int \e^{d(x,o)/D} \,\d\nu(x)\leq
(4+J(o)^2/s^2) \,\e^{m/D}
\]
and so
for any $1$-Lipschitz
function $f:X\to \R$ and $t\geq0$ we have
\[
\Pr\left(
\abs{f-f(o)}\geq t+m\right)\leq (8+2J(o)^2/s^2) \,\e^{-t/D}
\]
\end{thm}

So we get exponential concentration with caracteristic decay distance
$s^2/\rho$.

Note that the last assumption is satisfied with $s=2\sigma_\infty$ thanks to
Proposition 1.16 in~\cite{Led01}.

Before proceeding to the proof, let us show how this applies to the
geometric distribution above on $\N$. We take of course $o=0$ and $r=1$.
We can take $s=2\sigma_\infty=2$. Now there is only one point $x$ with $r\leq
d(o,x)<2r$, which is $x=1$. It satisfies $m_1=p\delta_0+(1-p)\delta_2$,
so that $\Td(m_1,\delta_0)=2(1-p)$, which is smaller than $d(0,1)=1$ if
and only if $p<1/2$ as was to be expected. So we can take
$\rho=1-2(1-p)=2p-1$. We get exponential concentration with
characteristic distance $4/(2p-1)$. When $p$ is very close to $1$ this is
not so good (because the discretization is too coarse), but when $p$ is
close to $1/2$ this is within a factor $2$ of the optimal value.

Another example is the stochastic differential equation $\d X_t=S\,
\d B_t-\alpha \frac{X_t}{\abs{X_t}}\,\d t$ on $\R^n$, for which $\exp
(-\abs{x}\alpha/S^2)$ is a reversible measure. Consider the Euler
scheme at time $\deltat$ for this stochastic differential equation.
Taking $r=nS^2/\alpha$ yields that $\rho\geq \alpha \,\deltat/2$
after some simple computation. Since we have $s^2=S^2 \deltat$ for
Gaussian measures at time $\deltat$, we get exponential concentration
with characteristic decay distance $2S^2/\alpha$, which is correct up
to a factor $2$. The additive constant in the deviation inequality is
$m=r+\rho(1+J(o)^2/4s^2)+2s^2/\rho$ which is equal to
$(n+4)S^2/\alpha+O(\deltat)$ (note that $J(o)^2\approx s^2$), which
is the correct order of magnitude for the average distance to $0$ in
dimension $n$.

If $\kappa>0$ in some large enough ball around $o$, then the invariant
distribution is unique. However, this is not true in general:
for example, start with the random walk on $\N$ above with a geometric
invariant distribution; now consider the disjoint union $\N\cup
(\N+\frac12)$ where we keep the same random walk on $\N$ and the same
walk translated by $\frac12$ on $\N+\frac12$: clearly there are two
disjoint invariant distributions, however, curvature is
non-negative and the assumptions of the theorem are satisfied with $r=1$
and $o=0$.

\begin{dem}[ of the theorem]

Let us first prove a lemma which shows how non-negative curvature
transmits the ``pulling''.

\begin{lem}
Let $x\in X$ with $d(x,o)\geq r$. Then $\Td(m_x,o)\leq d(x,o)-\rho$.
\end{lem}

\begin{dem}
If $d(o,x)<2r$ then this is one of the assumptions. So we suppose that
$d(o,x)\geq 2r$.

Since $X$ is $r$-geodesic, let $o=y_0,y_1,y_2,\ldots,y_n=x$ be a sequence
of points with $d(y_i,y_{i+1})\leq r$ and $\sum d(y_i,y_{i+1})=d(o,x)$.
We can assume that $d(o,y_2)>r$ (otherwise, remove $y_1$). 
Set $z=y_1$ if $d(o,y_1)=r$ and $z=y_2$ if $d(o,y_1)<r$, so that $r\leq
d(o,z)<2r$.
Now
\begin{align*}
\Td(\delta_o,m_x)
&\leq \Td(\delta_o,m_z)+\Td(m_z,m_x)
\\&\leq
d(o,z)-\rho+d(z,x)
\end{align*}
since $\kappa(z,x)\geq 0$. The conclusion
follows from the fact that $d(o,x)=d(o,z)+d(z,x)$.
\end{dem}

We are now ready to prove the theorem. The idea is to consider the
function $\e^{\lambda d(x,o)}$. For points far away from the origin, since
under the random walk the average distance to the origin decreases by
$\rho$ by the previous lemma, we expect the function to be multiplied by
$\e^{-\lambda \rho}$ under the random walk operator. Close to the origin,
the evolution of the function is controlled by the variance $s^2$ and the
jump $J(o)$ of the origin. Since the integral of the function is
preserved by the random walk operator, and it is multiplied by a quantity
$<1$ far away, this shows that the weight of faraway points cannot be too
large.

More precisely, we need to tamper a little bit with what happens around
the origin.
Let $\phi:\R_+\to\R_+$ be defined by $\phi(x)=0$ if $x<r$;
$\phi(x)=(x-r)^2/kr$ if $r\leq x< r(\frac{k}{2}+1)$ and $\phi(x)=x-r-kr/4$
if $x\geq r(\frac{k}{2}+1)$, for some $k>0$ to be chosen later. Note that
$\phi$ is a $1$-Lipschitz function and that $\phi''\leq 2/kr$.

If $Y$ is any random variable with values in $\R_+$, we have
\[
\E \phi(Y)\leq \phi(\E Y)+\frac12 \Var Y\,\sup \phi''\leq \phi(\E
Y)+\frac{1}{kr}\Var Y
\]

Now choose some $\lambda>0$ and consider the function $f:X\to \R$ defined
by $f(x)=\e^{\lambda \phi(d(o,x))}$. Note that $\phi(d(o,x))$ is
$1$-Lipschitz, so that by the Laplace transform assumption we have
\[
M\!f(x)\leq \e^{\lambda^2s^2/2}\e^{\lambda M\phi(d(o,x))}
\]

The Laplace transform assumption implies that the variance under $m_x$ of
any $1$-Lipschitz function is at most $s^2$. So by the remark above, we have
\[
M\phi(d(o,x))\leq \phi(\Td(m_x,\delta_o))+\frac{s^2}{kr}
\]
so that finally
\[
M\!f(x)\leq \e^{\lambda^2s^2/2+\lambda s^2/kr}\e^{\lambda \Td(m_x,\delta_o)}
\]

So for any $x$ with $d(o,x)\geq r$, we get
\[
M\!f(x)\leq \e^{\lambda^2s^2/2+\lambda s^2/kr}\e^{\lambda \phi(d(x,o)-\rho)}
\]

If $d(x,o)\geq r(\frac{k}{2}+1)+\rho$ then
$\phi(d(x,o)-\rho)=\phi(d(x,o))-\rho$ so that
\[
M\!f(x)\leq \e^{\lambda^2s^2/2+\lambda s^2/kr-\lambda \rho} f(x)
\]

If $r\leq d(x,o)<r(\frac{k}{2}+1)+\rho$, then $\phi(d(x,o)-\rho)\leq
\phi(d(x,o))$ so that
\[
M\!f(x)\leq \e^{\lambda^2s^2/2+\lambda s^2/kr} f(x)
\]

If, finally, $d(x,o)<r$, then use non-negative curvature 
to write $\Td(m_x,\delta_o)\leq
\Td(m_x,m_o)+J(o)\leq d(x,o)+J(o)$ so that $\phi(\Td(m_x,\delta_o))\leq
\phi(r+J(o))=J(o)^2/kr$ and
\[
M\!f(x)\leq \e^{\lambda^2s^2/2+\lambda s^2/kr+\lambda J(o)^2/kr}f(x)
\]

Let $\nu$ be a probability measure such that $\int f\,\d\nu<\infty$. 
Let $X'=\{x\in X,\; d(x,o)<r(\frac{k}{2}+1)\}$ and $X''=X\setminus X'$.
Set
$A(\nu)=\int_{X'} f\,\d\nu$ and
$B(\nu)=\int_{X''} f\,\d\nu$.
We have shown that
\begin{align*}
\int f \,\d(\nu\ast m)&=\int M\!f\,\d\nu=\int_{X'} M\!f\,\d\nu
+\int_{X''} M\!f\,\d\nu
\\&\leq \e^{\lambda^2s^2/2+\lambda s^2/kr+\lambda J(o)^2/kr} \int_{X'}f\,\d\nu
+ \e^{\lambda^2 s^2/2+\lambda s^2/kr-\lambda \rho} \int_{X''} f\,\d\nu
\end{align*}
so that
\[
A(\nu\ast m)+B(\nu\ast m)\leq \alpha A(\nu)+
\beta B(\nu)
\]
with $\alpha=\e^{\lambda^2s^2/2+\lambda s^2/kr+\lambda J(o)^2/kr}$ and $\beta=\e^{\lambda^2 s^2/2+\lambda s^2/kr-\lambda \rho}$.

Choose $\lambda$ small enough and $k$ large enough (see below) so that
$\beta<1$.
Using that $A(\nu)\leq
\e^{\lambda kr/4}$ for any measure $\nu$, we get $\alpha
A(\nu)+\beta B(\nu)\leq (\alpha-\beta)\e^{\lambda kr/4}+\beta
(A(\nu)+B(\nu))$. In particular, if $A(\nu)+B(\nu)\leq
\frac{(\alpha-\beta)\e^{\lambda r}}{1-\beta}$, we get $\alpha
A(\nu)+\beta B(\nu)\leq \frac{(\alpha-\beta)\e^{\lambda kr/4}}{1-\beta}$. So setting
$R=\frac{(\alpha-\beta)\e^{\lambda kr/4}}{1-\beta}$, we have just shown that
the set $C$ of probability measures $\nu$ such that $\int f\,\d\nu\leq R$ is
invariant under the random walk.

Moreover, if $A(\nu)+B(\nu)>R$ then $\alpha A(\nu)+\beta
B(\nu)<A(\nu)+B(\nu)$. Hence, if $\nu$ is a invariant distribution, necessarily $\nu\in C$. This, together with an evaluation of $R$
given below, provides the bound for $\int f\,\d\nu$ stated in the
theorem.

We now turn to existence of a invariant distribution.
First, $C$ is obviously closed and convex. Moreover, $C$ is tight: indeed
if $K$ is a compact, say included in a ball of radius $a$ around $o$,
then for any $\nu\in C$ we have $\nu(X\setminus K)\leq R \e^{-\lambda
a}$. So by Prokhorov's theorem, $C$ is compact in the weak convergence
topology. So $C$ is compact convex in the topological vector space of all
(signed) Borel measures on $X$, and is invariant by the random walk
operator, which is an affine map. By the Markov--Kakutani theorem
(Theorem~I.3.3.1 in~\cite{GD03}), it has a fixed point.

Let us finally evaluate $R$. We have
\begin{align*}
R&=\frac{\alpha/\beta-1}{1/\beta-1}\,\e^{\lambda kr/4}
\\&=\frac{\e^{\lambda J(o)^2/kr+\lambda\rho}-1}
{
\e^{\lambda \rho-\lambda s^2/kr-\lambda^2s^2/2}-1
}
\,\e^{\lambda kr/4}
\\&\leq 
\frac{\rho+J(o)^2/kr}{\rho-s^2/kr-\lambda s^2/2}
\,\e^{\lambda J(o)^2/kr+\lambda\rho+\lambda kr/4}
\end{align*}
using $\e^a-1\leq ae^a$ and $\e^a-1\geq a$.

Now take $\lambda=\rho/s^2$ and $k=4s^2/r\rho$. This yields
\[
R\leq (4+J(o)^2/s^2)\,\e^{\lambda(s^2/\rho+\rho(1+J(o)^2/4s^2))}
\]

Let $\nu$ be some invariant distribution. Since $d(x,o)\leq
\phi(d(x,o)+r(1+k/4))$ we have $\int \e^{\lambda d(x,o)} \,\d\nu\leq
\e^{\lambda r(1+k/4)}\int f\,\d\nu\leq Re^{\lambda r(1+k/4)}$ hence the
result in the theorem.
\end{dem}

\section{Ricci curvature and Gromov--Hausdorff topology}
\label{sec:GH}

We introduce here a Gromov--Hausdorff-like topology for metric spaces
equipped with a random walk. Two spaces are close in this topology if
they are close in the Gromov--Hausdorff topology and if moreover, the
measures issuing from each point $x$ are (uniformly) close in the $L^1$
transportation distance. More precisely:

\begin{defi}
Let $\left(X,(m_x)_{x\in X}\right)$ and $\left(Y,(m_y)_{y \in Y}\right)$ be two metric spaces
equipped with a random walk. For $e>0$, we say that these spaces are
$e$-close if there exists a metric space $Z$ and two isometric embeddings
$f_X:X\hookrightarrow Z$, $f_y:Y\hookrightarrow Z$ such that the Hausdorff distance between
$f_X(X)$ and $f_Y(Y)$ is at most $e$, and, moreover, for any $x\in X$,
there exists $y\in Y$ such that $d_Z(f_X(x),f_Y(y))\leq e$ and the $L^1$
transportation distance between the pushforward measures
$f_X(m_X)$ and $f_Y(m_y)$ is at most $2e$, and likewise for any $y\in Y$.
\end{defi}

The Ricci curvature is a continuous function in this topology. Namely, a
limit of spaces with Ricci curvature at least $\kappa$ has Ricci
curvature at least $\kappa$.

Below, we will relax the definition of Ricci curvature so as to allow any
variation at small scale; withthis perturbed definition, having Ricci
curvature greater than $\kappa$ will become an \emph{open} property. In
particular, any space close to a space with positive Ricci curvature will
have positive Ricci curvature in this perturbed sense.

\begin{prop}
Let $\left(X^n,(m^N_x)_{x\in X^N}\right)$ be a sequence of metric spaces
with random walk, converging to a metric space with random walk
$\left(X,(m_x)_{x\in X}\right)$. Let $x,y$ be two distinct points in $X$
and let $(x^N,y^N)\in (X^N,Y^N)$ be a sequence of pairs of points
converging to $(x,y)$.
Then $\kappa(x^N,y^N)\to \kappa(x,y)$.

In particular, if all spaces $X^N$ have Ricci curvature at least
$\kappa$, then so does $X$.
\end{prop}

In order for positive curvature to be an open property in some topology à
la Gromov--Hausdorff, one needs a rougher behavior at small scales. This is
achieved as follows.
 
\begin{defi}
Let $(X,d)$ be a metric space equipped with a random walk $m$. Let
$\delta\geq 0$. The \emph{Ricci curvature up to $\delta$ along $x,y\in
X$} is
\[
\kappa^\delta(x,y):=1-\frac{(\Td(m_x,m_y)-\delta)_+}{d(x,y)}
\]
i.e.~it is the largest $\kappa\leq 1$ for which one has
\[
\Td(m_x,m_y)\leq (1-\kappa)d(x,y)+\delta
\]
\end{defi}

With this definition, the following is easy.

\begin{prop}
Let $(X,(m_x))$ be a metric space with random walk with Ricci curvature
at least $\kappa$ up to $\delta\geq 0$. 
Let $\delta'>0$. Then there exists a neighborhood
$\mathcal{V}_X$ of $X$ such that any space $Y\in \mathcal{V}_X$ has Ricci
curvature at least $\kappa$ up to $\delta+\delta'$.
\end{prop}

Consequently, the property ``having curvature at least $\kappa$ for some
$\delta\geq 0$'' is open.

\section{$L^2$ Bonnet--Myers theorems}

\label{sec:strongBM}

As seen in Section~\ref{sec:weakBM}, it is generally not possible to give
a bound for the diameter of a positively curved space involving the
square root of curvature, because of such simple counterexamples as the
discrete cube.  Here we describe additional conditions which provide such
a bound in two different types of situation.

We first give a bound similar to the Bonnet--Myers one, but on the
\emph{average} distance between two points rather than the diameter; it
holds when there is an ``attractive point'' and is relevant for examples
such as the Ornstein--Uhlenbeck process (Example~\ref{ex:OU}) or its
discrete analogue (Example~\ref{ex:discrOU}).

Next, we give a direct generalization of the genuine Bonnet--Myers
theorem for Riemannian manifolds.  Actually, the only example where a
Bonnet--Myers theorem holds seems to be the ordinary Brownian motion on a
Riemannian manifold.
Despite this lack of further
examples, we found it interesting to provide an axiomatization of the
Bonnet--Myers theorem in our language.  This is done by reinforcing the
positive curvature assumption, which compares the transportation distance
between the measures issuing from two points $x$ and $y$ at a given time,
by requiring a transportation distance inequality between the measures
issuing from two given points at \emph{different} times.

\subsection{Average $L^2$ Bonnet--Myers}

We now describe a Bonnet--Myers-like estimate on the average distance
between two points, provided there is some ``attractive point''. This is
rather similar to Theorem~\ref{thm:expconc} in non-negative curvature.

\begin{prop}[ (Average $L^2$ Bonnet--Myers)]
Let $(X,d,(m_x))$ be a metric space with random walk, with Ricci
curvature at least $\kappa>0$. Suppose that for some $o\in X$ and $r\geq
0$, one has
\[
\int d(o,y)\,\d m_x(y)\leq d(o,x)
\]
for any $x\in X$ with $r\leq d(o,x)<2r$, and that moreover $X$ is
$r$-geodesic.

Then
\[
\int d(o,x)\,\d\nu(x)\leq \sqrt{\frac{1}{\kappa}\int
\frac{\sigma(x)^2}{n_x}\,\d\nu(x)}+5r
\]
where as usual $\nu$ is the invariant distribution.
\end{prop}

Note that the assumption $\int d(o,y)\,\d m_x(y)\leq d(o,x)$ cannot hold
for $x$ in some ball around $o$ unless $o$ is a fixed point. This is why
the assumption is restricted to an annulus.

As in the Gaussian concentration theorem (Theorem~\ref{thm:gaussconc}),
in case $\sigma(x)^2$ is Lipschitz, Corollary~\ref{cor:avlip} may provide
a useful bound on $\int \frac{\sigma(x)^2}{n_x}\,\d\nu(x)$ in terms of
its value at some point.

As a first example, consider the discrete Ornstein--Uhlenbeck process of
Example~\ref{ex:discrOU}, which is the Markov chain on $\{-N,\ldots,N\}$
given by the transition probabilities $p_{k,k}=1/2$, $p_{k,k+1}=1/4-k/4N$
and$p_{k,k-1}=1/4+k/4N$; the Ricci curvature is $\kappa=1/2N$, and the
invariant distribution is the binomial $\binom{2N}{N+k}$. This example is
interesting because the diameter is $2N$ (as is the bound provided by
Proposition~\ref{prop:weakBM}), whereas the average distance between two
points is $\approx \sqrt{N}$. It is immediate to check $0$ is attractive,
namely that $o=0$ and $r=1$ fulfill the assumptions. Since
$\sigma(x)^2\approx 1$ and $\kappa\approx 1/N$, the proposition recovers
the correct order of magnitude for distance to the origin.

Our next example is the Ornstein--Uhlenbeck process
$\d X_t=-\alpha \,X_t \,\d t+s\,\d B_t$ on $\R^N$ (Example~\ref{ex:OU}).
Here it is clear that $0$ is attractive in some sense, so $o=0$ is
a natural choice. The invariant distribution is a Gaussian of variance
$s^2/\alpha$; under this distribution the average distance to $0$ is
$\approx \sqrt{Ns^2/\alpha}$.

At small time $\tau$, a point $x\in \R^N$ is sent to a Gaussian centered
at $(1-\alpha\tau)x$, of variance $\tau s^2$. The average quadratic
distance to the origin under this Gaussian is
$(1-\alpha\tau)^2d(0,x)^2+Ns^2\tau+o(\tau)$ by a simple computation. If
$d(0,x)^2>Ns^2/2\alpha$ this is less than $d(0,x)^2$, so that we can take
$r=\sqrt{Ns^2/2\alpha}$. Considering the random walk discretized at time
$\tau$ we have we have $\kappa\sim \alpha \tau$, $\sigma(x)^2\sim
Ns^2\tau$ and $n_x\approx N$. So in the proposition above, the first term
is $\approx \sqrt{s^2/\alpha}$, whereas the second term is $5r\approx
\sqrt{Ns^2/\alpha}$, which is thus dominant. So the proposition gives the
correct order of magnitude; in this precise case, the first term in the
proposition reflects concentration of measure (which is
dimension-independent for Gaussians), whereas it is the second term $5r$
which carries the correct dependency on dimension for the average
distance to the origin.

\begin{dem}
Let $\phi:\R\to \R$ be the function defined by $\phi(x)=0$ if $x\leq 2r$,
and $\phi(x)=(x-2r)^2$ otherwise. Note that for any real-valued random
variable $Y$, we have
\[
\E \phi(Y)\leq \phi(\E Y)+\frac12 \Var Y \,\sup \phi''=\phi(\E Y)+\Var Y
\]

Now let $f:X\to \R$ be defined by $f(x)=\phi(d(o,x))$. We are going to
show that 
\[
M\!f(x)\leq (1-\kappa)^2f(x)+\frac{\sigma(x)^2}{n_x}+9r^2
\]
for all $x\in X$.
Since
$\int f\,\d\nu=\int M\!f\,\d\nu$, we will get $\int
f\,\d\nu\leq(1-\kappa)^2\int f\,\d\nu+\int
\frac{\sigma(x)^2}{n_x}\,\d\nu+9r^2$ which easily implies the result.

First, suppose that $r\leq d(o,x)<2r$. We have $f(x)=0$. Now
$\int d(o,y)\,\d m_x(y)$ is at most $d(o,y)$ by assumption. Using the
bound above for $\phi$, together with the definition of $\sigma(x)^2$ and
$n_x$, we get
\[
M\!f(x)=\int \phi(d(o,y))\,\d m_x(y)\leq \phi\left(\int d(o,y)\,\d
m_x(y)\right)+\frac{\sigma(x)^2}{n_x}=\frac{\sigma(x)^2}{n_x}
\]
since $\int d(o,y)\,\d m_x(y)\leq 2r$ by assumption.

Second, suppose that $d(x,o)\geq 2r$. Using that $X$ is $r$-geodesic, we
can find a point $x'$ such that $d(o,x)=d(o,x')+d(x',x)$ and $r\leq
d(o,x')<2r$ (take the second point in a sequence joining $o$ to $x$).
Now we have
\begin{align*}
\int  d(o,y)\,\d m_x(y)&=\Td(\delta_o,m_x)
\\&\leq \Td(\delta_o,m_{x'})+\Td(m_{x'},m_x)
\\&\leq \Td(\delta_o,m_{x'})+(1-\kappa) d(x',x)
\\&=\int  d(o,y)\,\d m_{x'}(y)+(1-\kappa) d(x',x)
\\&\leq d(o,x')+(1-\kappa)d(x',x)\leq (1-\kappa)d(o,x)+2\kappa r
\end{align*}
and as above, this implies
\begin{align*}
M\!f(x)&\leq \phi\left(\int  d(o,y)\,\d
m_x(y)\right)+\frac{\sigma(x)^2}{n_x}
\\&\leq \left((1-\kappa)d(o,x)+2\kappa
r-2r\right)^2+\frac{\sigma(x)^2}{n_x}
\\&=(1-\kappa)^2\phi(d(o,x))+\frac{\sigma(x)^2}{n_x}
\end{align*}
as needed.

The last case to consider is $d(o,x)<r$. In this case we have 
\begin{align*}
\int
d(o,y)\,\d m_x(y)&=\Td(\delta_o,m_x)
\\&\leq \Td(\delta_o,m_o)+\Td(m_o,m_x)=J(o)+\Td(m_o,m_x)
\\&\leq J(o)+(1-\kappa)d(o,x)\leq J(o)+r
\end{align*}

So we need to bound $J(o)$. If $X$ is included in the ball of radius
$r$ around $o$, the result trivially holds, so that we can assume that
there exists a point $x$ with $d(o,x)\geq r$. Since $X$ is $r$-geodesic
we can assume that $d(o,x)< 2r$ as well. Now
$J(o)=\Td(m_o,\delta_o)\leq \Td(m_o,m_x)+\Td(m_x,\delta_o)\leq
(1-\kappa)d(o,x)+\Td(m_x,\delta_o)\leq (1-\kappa)d(o,x)+d(o,x)$ by
assumption, so that $J(o)\leq 4r$.

Plugging this into the above, for $d(o,x)<r$ we get $\int
d(o,y)\,\d m_x(y)\leq 5r$ so that $\phi(\int
d(o,y)\,\d m_x(y))\leq 9r^2$ hence $M\!f(x)\leq
9r^2+\frac{\sigma(x)^2}{n_x}$.

Combining the results, we get that whatever $x\in X$
\[
M\!f(x)\leq (1-\kappa)^2f(x)+\frac{\sigma(x)^2}{n_x}+9r^2
\]
as needed.
\end{dem}

\subsection{Strong $L^2$ Bonnet--Myers}

As mentioned above, positive Ricci curvature alone does not imply a
$1/\sqrt{\kappa}$-like diameter control, because of such simple
counter-examples as the discrete cube or the Ornstein--Uhlenbeck process.
We now extract a property satisfied by the ordinary Brownian motion on
Riemannian manifolds (without drift), which guarantees a genuine
Bonnet--Myers theorem. Of course, this is of limited interest since the
only available example is Riemannian manifolds, but nevertheless we found
it interesting to find a sufficient condition expressed in our present
language.

Our definition of Ricci curvature controls the transportation distance
between the measures issuing from two points $x$ and $x'$ at a given time
$t$. The condition we will now use controls the transportation distance
between the measures issuing from two points at two \emph{different} times. It
is based on what holds for Gaussian measures in $\R^N$. For any $x,x' \in
\R^N$ and $t,t'>0$, let $m_x^{\ast t}$ and $m_{x'}^{\ast t'}$ be the laws
of the standard Brownian motion issuing from $x$ at time $t$ and from $x'$ at time
$t'$, respectively. It is easy to check that the $L^2$ transportation
distance between these two measures is
\[
\Tdd(m_x^{\ast t},m_{x'}^{\ast
t'})^2=d(x,x')^2+N(\sqrt{t}-\sqrt{t'})^2
\]
hence
\[
\Td(m_x^{\ast t},m_{x'}^{\ast
t'})\leq d(x,x')+\frac{N(\sqrt{t}-\sqrt{t'})^2}{2d(x,x')}
\]

The important feature here is that, when $t'$ tends to $t$, the second
term is of second order in $t'-t$. This is no more the case if we
add a drift term to the diffusion.

We now take this inequality as an assumption and use it to mimick the
traditional proof of the Bonnet--Myers theorem. Here, for simplicity of
notation we suppose that we are given a continuous-time Markov chain;
however, the proof uses only a finite number of different values of $t$,
so that discretization is possible (this is important in Riemannian
manifolds, because the heat kernel is positive on the whole manifold at
any positive time, and there is no simple control on it far away from the
initial point; taking a discrete approximation with bounded steps solves
this problem).

\begin{prop}[ (Strong $L^2$ Bonnet--Myers)]
Let $X$ be a metric space equipped with a continuous-time random walk
$m^{\ast t}$. Assume that $X$ is $\eps$-geodesic, and that there exists
constants $\kappa>0, C\geq 0$ such that for any two
small enough $t,t'$, for any
$x,x'\in X$ with $\eps\leq d(x,x')\leq 2\eps$ one has
\[
\Td(m_x^{\ast t},m_{x'}^{\ast
t'})\leq \e^{-\kappa\inf (t,t')}d(x,x')+\frac{C(\sqrt{t}-\sqrt{t'})^2}{2d(x,x')}
\]
with $\kappa>0$. Assume moreover that $\eps\leq \frac12
\sqrt{C/2\kappa}$.

Then
\[
\diam X\leq \pi
\sqrt{\frac{C}{2\kappa}}\,\left(1+\frac{4\eps}{\sqrt{C/2\kappa}}\right)
\]
\end{prop}

When $t=t'$, the assumption reduces to $\Td(m_x^{\ast t},m_{x'}^{\ast
t})\leq \e^{-\kappa t}d(x,x')$, which is just the continuous-time version
of the positive curvature assumption.
The constant $C$ plays the role of
a diffusion constant, and is equal to $N$ for (a discrete approximation
of) Brownian motion on a Riemannian manifold. We restrict the assumption
to $d(x,x')\geq \eps$ to avoid divergence problems for
$\frac{C(\sqrt{t}-\sqrt{t'})^2}{2d(x,x')}$ when $x'\to x$.

For the Brownian motion on an $N$-dimensional Riemannian manifold, we can take
$\kappa=\frac12\inf \Ric$ by Bakry-Émery theory (the $\frac12$ is due to
the fact that the infinitesimal generator of Brownian motion is
$\frac12\Delta$), and $C=N$ as in $\R^N$. So we get the usual
Bonnet--Myers theorem, up to a factor $\sqrt{N}$ instead of $\sqrt{N-1}$
(similarly to our spectral gap estimate in comparison with the
Lichnerowicz theorem),
but with the correct constant $\pi$.

\begin{dem}
Let $x,x'\in X$. Since $X$ is $\eps$-geodesic, we can find a sequence
$x=x_0,x_1,\ldots,x_{k-1},x_k=x'$ of points in $X$ with
$d(x_i,x_{i+1})\leq \eps$ and $\sum d(x_i,x_{i+1})=d(x_0,x_k)$.
By taking
a subsequence (denoted $x_i$ again), we can assume that $\eps\leq
d(x_i,x_{i+1})\leq 2\eps$ instead. 

Set $t_i=\eta\sin\left(\frac{\pi d(x,x_i)}{d(x,x')}\right)^2$ for
some (small) value of $\eta$ to be chosen later. Now, since $t_0=t_k=0$
we have
\begin{align*}
d(x,x')&=\Td(\delta_x,\delta_{x'})
\leq \sum \Td (m_{x_i}^{\ast t_i},m_{x_{i+1}}^{\ast t_{i+1}})
\\&\leq
\sum\,
\e^{-\kappa
\inf(t_i,t_{i+1})}d(x_i,x_{i+1})+\frac{C(\sqrt{t_{i+1}}-\sqrt{t_i})^2}{2d(x_i,x_{i+1})}
\end{align*}
by assumption. Now, for $a<b$ we have $\sin b-\sin a=2\sin \frac{b-a}{2}
\cos\frac{a+b}{2}\leq (b-a)\cos\frac{a+b}{2}$ so that
\[
\frac{C(\sqrt{t_{i+1}}-\sqrt{t_i})^2}{2d(x_i,x_{i+1})}\leq 
\frac{C\eta\pi^2 d(x_i,x_{i+1})}{2d(x,x')^2}\,\cos^2 \left(\pi\frac{
d(x,x_i)+d(x,x_{i+1})}{2d(x,x')}\right)
\]

Besides, if $\eta$ is small enough, one has $\e^{-\kappa
\inf(t_i,t_{i+1})}=1-\kappa \inf(t_i,t_{i+1})+O(\eta^2)$. So we get
\begin{align*}
d(x,x')\leq \sum & \,d(x_i,x_{i+1})-\kappa \inf(t_i,t_{i+1}) d(x_i,x_{i+1})
\\&+
\frac{C\eta\pi^2 d(x_i,x_{i+1})}{2d(x,x')^2}\,\cos^2 \left(\pi\frac{
d(x,x_i)+d(x,x_{i+1})}{2d(x,x')}\right)
+O(\eta^2)
\end{align*}

Now the terms $\sum
d(x_i,x_{i+1})\,\cos^2 \left(\pi\frac{
d(x,x_i)+d(x,x_{i+1})}{2d(x,x')}\right)$ and $\sum \inf(t_i,t_{i+1})
d(x_i,x_{i+1})$ are close to the integrals $d(x,x')\int_0^1 \cos^2(\pi u)\,\d
u$ and $d(x,x')\eta\int_0^1 \sin^2(\pi u)\,\d u$ respectively; the relative
error in the Riemann sum is easily bounded by $\pi\eps/d(x,x')$ so that
\begin{align*}
d(x,x')\leq & \,d(x,x')-\kappa \,\eta\,
d(x,x')\left(\frac12-\frac{\pi\eps}{d(x,x')}\right)
\\&+\frac{C\eta\pi^2}{2d(x,x')^2}\,d(x,x')\left(\frac12+\frac{\pi\eps}{d(x,x')}\right)
+O(\eta^2)
\end{align*}
hence, taking $\eta$ small enough,
\[
d(x,x')^2\leq
\frac{C\pi^2}{2\kappa}\,\frac{1+2\pi\eps/d(x,x')}{1-2\pi\eps/d(x,x')}
\]
so that either $d(x,x')\leq \pi\sqrt{C/2\kappa}$, or
$2\pi\eps/d(x,x')\leq 2\pi\eps/\pi\sqrt{C/2\kappa}\leq 1/2$ by the
assumption that $\eps$ is small, in which
case we use $(1+a)/(1-a)\leq 1+4a$ for $a\leq 1/2$, hence the conclusion.
\end{dem}

%
%
%
%

\section{Transportation distance in Riemannian manifolds}
\label{sec:riemproof}

Here we give the proofs of Proposition~\ref{prop:riemdist} and of the
statements of Example~\ref{ex:riem} and Section~\ref{sec:riem}.

We begin with Proposition~\ref{prop:riemdist} and evaluation of the Ricci
curvature of the random walk at scale $\eps$.

Let $X$ be a smooth $N$-dimensional Riemannian manifold and let $x\in X$.
Let $v,w$ be unit tangent vectors at $x$. Let $\delta,\eps>0$ small
enough. Let $y=\exp_x(\delta v)$. Let $x'=\exp_x (\eps w)$ and $y'=\exp_y
(\eps w')$ where $w'$ is the tangent vector at $y$ obtained by parallel
transport of $w$ along the geodesic $t\mapsto \exp_x (tv)$. The first
claim is that $d(x',y')=\delta\left(1-\frac{\eps^2}{2} K(v,w)+O(\delta
\eps^2 +\eps^3)\right)$.

We suppose for simplicity that $w$ and $w'$ are orthogonal to $v$.

We will work in cylindrical coordinates along the geodesic $t\mapsto
\exp_x(tv)$. Let $v_t=\frac{\d}{\d t} \exp_x(tv)$ be the speed of this
geodesic. Let $E_t$ be the orthogonal of $v_t$ in the tangent space at
$\exp_x(tv)$. Each point $z$ in some neighborhood of $x$ can be uniquely
written as $\exp_{\exp_x (\tau(z) v)} (\eps \zeta(z))$ for some $\tau(z) \in \R$ and
$\zeta(z)\in E_{\tau(z)}$.

Consider the function $f$ equal to the distance of a point to $\exp_x
(E_0)$ (taken in some small enough neighborhood of $x$), equipped with a
$-$ sign if the point is not on the same side of $E_0$ as $y$.  Clearly
$f$ is $1$-Lipschitz, so that $d(x',y')\geq f(y')-f(x')$.

The distance from $\exp_x(E_0)$ to $y'$ is realized by some geodesic
$\gamma$ starting at some point of $\exp_x(E_0)$ and ending at $y$. If
$\delta$ and $\eps$ are small enough, this geodesic is arbitrarily close
to the Euclidean situation so that the coordinate $\tau$ is strictly
increasing along $\gamma$. Let us parametrize $\gamma$ using the
coordinate $\tau$, so that $\tau(\gamma(t))=t$. Let also
$w_t=\zeta(\gamma(t))\in E_t$. In particular, $w_\delta=w'$.

Now by definition we have $\gamma(t)=\exp_{\exp_x (t v)} (\eps
w_t)$. Considering the family of geodesics $s\mapsto \exp_{\exp_x (t v)}
(s w_t)$ and applying the Jacobi equation yields
\[
\abs{\frac{\d \gamma(t)}{\d t}}^2=\abs{v_t}^2+2\eps \langle
v_t, \dot{w_t}\rangle+\eps^2 \abs{\dot{w_t}}^2-\eps^2\,\langle R(w_t,v_t)w_t,v_t\rangle +O(\eps^3)
\]
where $\dot{w_t}=\frac{D}{\d t} w_t$. But since by definition $w_t\in
E_t$, we have $\langle v_t, \dot{w_t}\rangle=0$. Since moreover
$\abs{v_t}=1$ we get
\[
\abs{\frac{\d \gamma(t)}{\d t}}=1+\frac{\eps^2}{2}
\abs{\dot{w_t}}^2-\frac{\eps^2}{2} \langle R(w_t,v_t)w_t,v_t\rangle
+O(\eps^3)
\]
which is always greater than $1-\frac{\eps^2}{2} \langle
R(w_t,v_t)w_t,v_t\rangle
+O(\eps^3)$. Integrating from $t=0$ to $t=\delta$ and using that $\langle
R(w_t,v_t)w_t,v_t\rangle= K(w,v) +O(\delta)$ yields that
the length of the geodesic is
\[
\delta\,\left(1-\frac{\eps^2}{2}\,K(v,w)+O(\eps^3)+O(\eps^2\delta)\right)
\]
so that the distance from $x'$ to $y'$ is at least this quantity. But
this value is achieved for $\dot{w_t}=0$, in which case $\gamma(0)=x'$ by
definition, so this is exactly $d(x',y')$. This proves
Proposition~\ref{prop:riemdist}.

\bigskip

Let us now prove the statement of Example~\ref{ex:riem}. Let $\mu_0,
\mu_1$ be the uniform probability measures on the balls of radius $\eps$
centered at $x$ and $y$ respectively. We have to prove that
\[
\Td(\mu_0,\mu_1)=d(x,y)\left(1-\frac{\eps^2}{2(N+2)}\,\Ric(v,v)\right)
\]
up to higher-order terms.

Let $\mu'_0,\mu'_1$ be the images under the exponential map, of the
uniform probability measures on the balls of radius $\eps$ in the tangent
spaces at $x$ and $y'$ respectively. So $\mu'_0$ is a measure having
density $1+O(\eps^2)$ w.r.t.\ $\mu_0$, and likewise for $\mu_1'$.

If we average Proposition~\ref{prop:riemdist} over $w$ in the ball of
radius $\eps$ in the tangent space at $x$, we get that
\[
\Td(\mu'_0,\mu'_1)\leq d(x,y)\left(1-\frac{\eps^2}{2(N+2)}\,\Ric(v,v)\right)
\]
up to higher-order terms,
since the coupling by parallel transport realizes this value. Indeed,
$\Ric(v,v)$ is the sum of $K(v,w)$ for $w$ in an orthonormal basis of
the tangent space at $x$. Consequently, the average of $K(v,w)$ on the
unit sphere is $\frac1N \Ric(v,v)$. Averaging on the ball instead of the
sphere yields an $\frac{1}{N+2}$ factor instead.

Now the density of $\mu'_0$, $\mu'_1$ with respect to $\mu_0$, $\mu_1$ is
$1+O(\eps^2)$. Moreover the $O(\eps^2)$ terms decompose as the sum of an
$O(d(x,y)\eps^2)$ term and an $O(\eps^2)$ term which is the same for
$\mu'_0$ and $\mu'_1$ (indeed, $\mu'_0$ and $\mu'_1$ coincide when $x=y$).
Plugging this in the estimate above, we get the inequality for
$\Td(\mu_0,\mu_1)$ up to higher-order terms.

The converse inequality is proven as follows: if $f$ is any $1$-Lipschitz
function, the $L^1$ transportation distance between measures $\mu_0$ and
$\mu_1$ is at least the difference of the integrals of $f$ under $\mu_0$
and $\mu_1$ (and actually, a clever choice of $f$ realizes this
transportation distance, see Theorem~1.14 in~\cite{Vil03}). Arguments
similar to the above for integrating under $\mu_0$ and $\mu_1$, applied
to the function $f$ above equal to the distance of a point to the set
$\exp_x(E_0)$, yield the desired inequality.

Finally, let us briefly sketch the proofs of the other statements of
Section~\ref{sec:riem}, namely, evaluation of the spread and local
dimension (Definition~\ref{def:jsd}). Up to a multiplicative factor $O(1+\eps)$, these can be
computed in the Euclidean space.

A simple computation shows that the expectation of the square distance of
two points taken at random in a ball of radius $\eps$ is
$\eps^2\,\frac{2N}{N+2}$, hence the value $\eps^2\,\frac{N}{N+2}$ for the
spread.

To evaluate the local dimension (Definition~\ref{def:jsd}), we have to
bound the maximal variance of a $1$-Lipschitz function on a ball of
radius $\eps$. We will prove that the local dimension $n_x$ is comprised
between $N-1$ and $N$. A projection to a coordinate axis provides a function with
variance $\frac{\eps^2}{N+2}$, so that local dimension is at most $N$.
For the other bound, let $f$ be a $1$-Lipschitz function on the ball and
let us compute an upper bound for its variance. Take $\eps=1$ for
simplicity. Write the ball of radius $1$ as the union of the spheres
$S_r$ of radii $r\leq 1$. Let $v(r)$ be the variance of $f$ restricted
to the sphere $S_r$, and let $a(r)$ be the average of $f$ on $S_r$. Then
associativity of variances gives
\[
\Var f=\int_{r=0}^1 v(r)\,\d\mu(r) +\Var_\mu a(r)
\]
where $\mu$ is the measure on the interval $[0;1]$ given by
$\frac{r^{N-1}}{Z}\,\d r$ with
$Z=\int_{r=0}^1 r^{N-1}\,\d r=\frac1N$.

Since the variance of a $1$-Lipschitz function on the $(N-1)$-dimensional
unit sphere is at most $\frac1N$, we have $v(r)\leq \frac{r^2}{N}$ so
that $\int_{r=0}^1 v(r)\,\d\mu(r)\leq \frac1{N+2}$. To evaluate the
second term, note that $a(r)$ is again $1$-Lipschitz as a function of
$r$, so that $\Var_\mu a(r)=\frac12 \iint
(a(r)-a(r'))^2\,\d\mu(r)\d\mu(r')$ is at most $\frac12 \iint
(r-r')^2\,\d\mu(r)\d\mu(r')=\frac{N}{(N+1)^2(N+2)}$. So finally
\[
\Var f\leq \frac1{N+2}+\frac{N}{(N+1)^2(N+2)}
\]
so that the local dimension $n_x$ is bounded below by
$\frac{N(N+1)^2}{N^2+3N+1}\geq N-1$.


\begin{thebibliography}{MMMM}

\bibitem[ABCFGMRS00]{LSI} C.~Ané, S.~Blachère, D.~Chafaï, P.~Fougères,
I.~Gentil, F.~Malrieu, C.~Roberto, G.~Scheffer, \emph{Sur les inégalités
de Sobolev logarithmiques}, Panoramas et Synthèses~\textbf{10}, Société
Mathématique de France (2000).

\bibitem[AMTU01]{AMTU01} A.~Arnold, P.~Markowich, G.~Toscani,
A.~Unterreiter, \emph{On convex Sobolev inequalities and the rate of
convergence to equilibrium for Fokker-Planck type equations},  Comm.
Partial Differential Equations~\textbf{26}  (2001),  n°~1-2, 43--100.

\bibitem[Ber03]{Ber03} M.~Berger, \emph{A panoramic view of Riemannian
geometry}, Springer, Berlin (2003).

\bibitem[BE84]{BE84} D.~Bakry, M.~Émery, \emph{Hypercontractivité de
semi-groupes de diffusion},  C.\ R.\ Acad.\ Sci.\ Paris Sér.\ I
Math.~\textbf{299} (1984), n°~15, 775--778.

\bibitem[BE85]{BE85} D.~Bakry, M.~Émery, \emph{Diffusions
hypercontractives},
Séminaire de probabilités, XIX, 1983/84.
Lecture Notes in Math.~\textbf{1123}, Springer, Berlin (1985),
177--206.

\bibitem[BL98]{BL98} S.~Bobkov, M.~Ledoux, \emph{On modified logarithmic
Sobolev inequalities for Bernoulli and Poisson measures},
J.\ Funct.\ Anal.~\textbf{156} (1998), n°~2, 347--365.

\bibitem[Bré99]{Bre99} P.~Brémaud, \emph{Markov chains}, Texts in Applied
Mathematics~\textbf{31}, Springer, New York (1999).

\bibitem[Che98]{Che98} M.-F.~Chen, \emph{Trilogy of couplings and general
formulas for lower bound of spectral gap}, in \emph{Probability towards
2000 (New York, 1995)}, Lecture Notes in Statist.~\textbf{128}, Springer,
New York (1998), 123--136.

\bibitem[CW97]{CW97} M.-F.~Chen, F.-Y.~Wang, \emph{Estimation of spectral
gap for elliptic operators}, Trans.\ Amer.\ Math.\ Soc.~\textbf{349}
(1997), n°~3, 1239--1267.

\bibitem[Dob56]{Dob56} R. L.~Dobru\v{s}in,   \emph{On the condition of the
central limit theorem for inhomogeneous Markov chains} (Russian), Dokl.\
Akad.\ Nauk SSSR (N.S.)~\textbf{108} (1956), 1004--1006. 

\bibitem[DGW04]{DGW04} H.~Djellout, A.~Guillin, L.~Wu,
\emph{Transportation cost-information inequalities and applications to
random dynamical systems and diffusions}, Ann.\ Prob.~\textbf{32} (2004),
n°~3B, 2702--2732.

\bibitem[GD03]{GD03} A.~Granas, J.~Dugundji, \emph{Fixed point theory},
Springer Monographs in Mathematics, Springer, New York (2003).

\bibitem[GH90]{GH90} É.~Ghys, P.~de la Harpe, \emph{Sur les groupes
hyperboliques d'après Mikhael Gromov}, Progress in Math.~\textbf{83},
Birkhäuser (1990).

\bibitem[Gri67]{Gri67} R. B.~Griffiths, \emph{Correlations in Ising
ferromagnets III}, Commun.\ Math.\ Phys.~\textbf{6} (1967), 121--127.

\bibitem[Gro86]{Gro86} M.~Gromov, in V.~Milman, G.~Schechtman,
\emph{Asymptotic theory of finite dimensional normed spaces}, Lecture
Notes in Mathematics~\textbf{1200}, Springer, Berlin (1986).

\bibitem[Jou]{Jou} A.~Joulin, \emph{Poisson-type deviation inequalities
for curved continuous time Markov chains}, preprint.

\bibitem[Led01]{Led01} M.~Ledoux, \emph{The concentration of measure
phenomenon}, Mathematical Surveys and Monographs~\textbf{89}, AMS (2001).

\bibitem[Lott]{Lott} J.~Lott, \emph{Optimal transport and Ricci curvature
for metric-measure spaces}, expository manuscript.

\bibitem[LV]{LV} J.~Lott, C.~Villani, \emph{Ricci curvature for
metric-measure spaces via optimal transport}, preprint.

\bibitem[Oht]{Oht} S.-i.~Ohta, \emph{On the measure contraction
property of metric measure spaces}, preprint.

\bibitem[Oli]{Oli} R.~I.~Oliveira, \emph{On the convergence to
equilibrium of Kac's random walk on matrices}, preprint,
\texttt{arXiv:0705.2253}

\bibitem[RS05]{RS05} M.-K.~von Renesse, K.-T.~Sturm, \emph{Transport
inequalities, gradient estimates, and Ricci curvature},
Comm.\ Pure Appl.\ Math.~\textbf{68} (2005), 923--940.

\bibitem[Sam]{Sam} M.~D.~Sammer, \emph{Aspects of mass transportation in
discrete concentration inequalities}, PhD thesis, Georgia institute of
technology, 2005,
\url{etd.gatech.edu/theses/available/etd-04112005-163457/unrestricted/sammer_marcus_d_200505_phd.pdf}

\bibitem[Stu06]{Stu06} K.-T.~Sturm, \emph{On the geometry of metric
measure spaces},  Acta Math.~\textbf{196} (2006),  n°1, 65--177.

\bibitem[Vil03]{Vil03} C.~Villani, \emph{Topics in optimal
transportation}, Graduate Studies in Mathematics~\textbf{58}, AMS
(2003).

\bibitem[Vil]{Vil} C.~Villani, \emph{Optimal transport, old and new},
July 12, 2007 version,
\url{www.umpa.ens-lyon.fr/~cvillani/Cedrif/B07B.StFlour.pdf}

\end{thebibliography}
\end{document}